\documentclass[11pt]{article}
\usepackage{amssymb,amsmath,amsthm,epsfig,times}

\def\C{{\mathbb C}}

\def\R{{\mathbb R}}
\def\*S{{\mathbb S}}

\def\Z{{\mathbb Z}}

\def\interior
{\begin{picture}(6,6)
\put (0,0){\line(1,0){6}}
\put (6,0){\line(0,1){6}}
\end{picture} \hspace{0.05in}}
\newcommand{\dis}{\displaystyle}

\newtheorem{theorem}{Theorem}
\newtheorem{proposition}{Proposition}[section]
\newtheorem{lemma}[proposition]{Lemma}

\newtheorem{definition}[proposition]{Definition}

\begin{document}

\large

\title{Residual $\bar\partial$-cohomology and the complex Radon transform on subvarieties of $\C P^n$.}
\author{Gennadi M. Henkin\footnote{Institut de Mathematiques, Universite Pierre et Marie Curie, 75252, BC247, Paris,
Cedex 05, France, and CEMI, Acad. Sc., 117418, Moscow, Russia, henkin@math.jussieu.fr}\ ,
Peter L. Polyakov\footnote{Department of Mathematics, University of Wyoming, Laramie, WY 82071, USA, polyakov@uwyo.edu}}
\maketitle

\begin{abstract}
We show that the complex Radon transform realizes an isomorphism between the quotient-space of
residual $\bar\partial$-cohomologies of a locally complete intersection algebraic subvariety in a linearly
concave domain of ${\C}P^n$ and the space of holomorphic solutions of the associated 
homogeneous system of differential equations with constant coefficients in the dual domain in $({\C}P^n)^*$.
\end{abstract}

\section{Introduction.}\label{Introduction}
\indent

In this article we consider two related problems: the first one is the description of
infinite-dimensional spaces of $\bar\partial$-cohomologies of subvarieties in
linearly concave domains of ${\C}P^n$ in terms of inverse Radon transform of the spaces of holomorphic
solutions of associated systems of differential equations in dual domains, and the second one is
the realization of the spaces of holomorphic
solutions of systems of linear differential equations in convex domains by Radon transforms 
of $\bar\partial$-cohomologies of associated subvarieties in dual domains.\\
\indent
The study of these problems was started by Martineau in \cite{Mar1}, \cite{Mar2}
and was continued in the papers \cite{GH}, \cite{HP1}, \cite{B}, \cite{He1}, \cite{He2}, \cite{DSc}.
The main result of Martineau in \cite{Mar1} was interpreted in \cite{GH} as the existence
of an isomorphism defined by the complex Radon transform between the space of
$(n,n-1)$ $\bar\partial$-cohomologies of a linearly concave domain $D\subset {\C}P^n$
and the space of holomorphic functions on the dual linearly convex domain
$D^*\subset ({\C}P^n)^*$.\\
\indent
We begin by describing the result that was produced by the study of the problems mentioned above in
\cite{HP1},\cite{He1} for the case of complex submanifolds in linearly concave domains in $\C P^n$.\\
\indent
Let $(z_0,\ldots,z_n)$ and $(\xi_0,\ldots,\xi_n)$ be the homogeneous
coordinates of points $z\in {\C}P^n$ and $\xi\in ({\C}P^n)^*$. Let
$\langle\xi\cdot z\rangle\stackrel{\rm def}{=}\sum\limits_{k=0}^n\xi_k\cdot z_k$, and let  $\C P^{n-1}_{\xi}$
denote the hyperplane
\begin{equation*}
\C P^{n-1}_{\xi}=\{z\in \C P^n:\ \langle\xi\cdot z\rangle=0\}.
\end{equation*}
\indent
Following  \cite{Mar1} and \cite{GH} we call a domain $D\subset \C P^n$ a linearly concave domain,
if there exists a continuous family of hyperplanes  $\C P^{n-1}(z)\subset D$ defined for $z\in D$ and satisfying
$z\in \C P^{n-1}(z)$.
We notice that in the original definition of linearly concave domains in \cite{Mar1}
the continuity of the family was not required, but the main results of \cite{Mar1}, \cite{Mar2},
\cite{GH}, \cite{HP1}, \cite{He1} are valid only under the assumption of existence of such family.\\
\indent
The following theorem was obtained in \cite{He1}.
%********************************************************************************
%******* Theorem OldTheorem  ******************************************************
\begin{theorem}\label{OldTheorem}
Let $D$ be a linearly concave domain in $\C P^n$, $n\ge 2$, and let $D^*\subset  ({\C}P^n)^*$be the dual domain
$$D^*=\left\{\xi\in  ({\C}P^n)^*: \C P^{n-1}_{\xi}\subset D\right\}.$$
Let $V$ be a $(n-m)$-dimensional connected algebraic manifold of the form
\begin{equation*}
V=\left\{z\in \C P^n:\ P_1(z)=\ldots =P_r(z)=0\right\},
\end{equation*}
where homogeneous polynomials $P_1,\ldots,P_r$ are such that everywhere on $V$
\begin{equation*}
\text{rank}\,\left[\text{grad}\,P_1,\ldots,\text{grad}\,P_r\right]=m.
\end{equation*}
Let $V_D=V\cap D$, let $Z^{(n-m,n-m-1)}(V_D)$ denote the space of
$\bar\partial$-closed smooth forms on $V_D$ of bidegree $(n-m,n-m-1)$,
and let $H^0\left(D^*\right)$ and $H^{(1,0)}\left(D^*\right)$ denote the spaces of holomorphic functions
and respectively holomorphic 1-forms on $D^*$.\\
\indent
Then the Radon transform
$${\cal R}_V:Z^{(n-m,n-m-1)}\left(V_D\right)\to H^{(1,0)}\left(D^*\right)$$
defined by the formula
\begin{equation}\label{OldRadon}
{\cal R}_V[\phi](\xi)=\sum_{j=0}^n\left(\int\limits_{z\in \C P^{n-1}_{\xi}\cap V}\langle\xi\cdot dz\rangle
\interior z_j \phi\right)d\xi_j
\end{equation}
induces a continuous linear operator on the space of cohomologies
$${\cal R}_V:H^{(n-m,n-m-1)}\left(V_D\right)\to H^{(1,0)}\left(D^*\right).$$
\indent
The following properties are satisfied:
\begin{itemize}
\item[(i)]
the subspace ${Ker}{\cal R}_V\subset H^{(n-m,n-m-1)}(V_D)$ is finite-dimensional  and consists of
restrictions to $V_D$ of $\bar\partial$-cohomologies from $H^{(n-m,n-m-1)}(V)$,
\item[(ii)]
the image of ${\cal R}_V$ is the following subspace in $H^{(1,0)}\left(D^*\right)$
\begin{multline}\label{OldImage}
{\cal R}_V\left(H^{(n-m,n-m-1)}\left(V_D\right)\right)\\
=\left\{f\in H^{(1,0)}\left(D^*\right):
f=dg\ \mbox{with}\ g\in H^0\left(D^*\right)\ \mbox{such that}\
\left\{P_k\left(\frac{\partial}{\partial\xi}\right)g=0\right\}_1^r\right\}.
\end{multline}
\end{itemize}
\end{theorem}
%********************************************************************************
{\bf Remarks.}
\begin{itemize}
\item
If $V\subset \C P^n$ is a smooth complete intersection, and $D\subset \C P^n$ is a linearly concave domain,
then in (\cite{He1} Theorem 5.1) an explicit inversion formula for $R_{V}$ is obtained in the spirit
of explicit fundamental principle of \cite{BP}.
\item
For $m=n-1$ the statement (i) of Theorem~\ref{OldTheorem} is a corollary of the inverse Abel theorem
(see Saint-Donat \cite{SD}, Griffiths \cite{Gr}). For $m<n-1$ and $V$ - complete intersection, the statement (i)
of Theorem~\ref{OldTheorem} is a consequence of Theorem 3.3 from \cite{HP1}.
\item
In the statement (ii) of Theorem~\ref{OldTheorem} if $\phi\in Z^{(n-m,n-m-1)}\left(V_D\right)$ is such that
$\phi=\partial\psi$ for $\psi\in Z^{(n-m-1,n-m-1)}\left(V_D\right)$, then $g$ is the image of $\psi$ under
the map introduced by Andreotti and Norguet (see \cite{AN}, \cite{O}).
\end{itemize}

\indent
The main result of this article is a natural generalization of Theorem~\ref{OldTheorem} to the
case of an arbitrary locally complete intersection in a linearly concave domain. In order to formulate
this theorem we need to introduce some additional definitions and notations.\\
\indent
Throughout the whole article we will denote by $D\subset \C P^n$ a linearly concave domain and by
$G=\C P^n\setminus D$ its complement. We will also denote by $D_{\delta}$ linearly concave subdomains
of $D$ with smooth boundaries $bD_{\delta}$ such that
\begin{equation*}
D_{\delta}\subset D_{\nu}\quad \text{for}\ \nu<\delta,\quad
\text{and}\quad \bigcup_{\delta} D_{\delta}=D.
\end{equation*}
The existence of a sequence of subdomains with the above properties is proved in Proposition~\ref{Exhaustion}.
We will denote by $G_{\delta}=\C P^n\setminus D_{\delta}\supset G$ and by $\mathring{G}=G\setminus bG$.

%***************************************************************************************
%******* Definition  LocallyComplete ***************************************************
\begin{definition}\label{LocallyComplete} (Locally Complete Intersections)
An analytic subvariety $V\subset \C P^n$ is called a locally complete intersection subvariety
in $\C P^n$ of pure dimension $n-m$ if there exist a finite open cover
$\left\{U_{\alpha}\right\}_{\alpha=1}^N$ of $\C P^n$  and collections of holomorphic functions
$\left\{F^{(\alpha)}_k\right\}$ in $U_{\alpha}$, such that
\begin{equation}\label{V_cap_U}
V\cap U_{\alpha}=\left\{z\in U_{\alpha}:\ F^{(\alpha)}_1(z)=\cdots=F^{(\alpha)}_m(z)=0\right\}
\end{equation}
with the structure sheaf ${\cal O}/{\cal I}$, where ${\cal O}$ is the structure sheaf of $\C P^n$, and
${\cal I}$ is the sheaf of ideals defined by polynomials $\{F^{\alpha}_k\}_{k=1}^m$.
\end{definition}
%***************************************************************************************
%***************************************************************************************
\indent
In our construction of  $\bar\partial$-closed residual currents on a locally complete intersection variety $V$
we will use a special vector bundle, the so-called {\it conormal vector bundle}.
To describe this bundle we consider a domain $U\subset \C P^n$,
a finite cover $\left\{U_{\alpha}\right\}_{\alpha=1}^N$ of $U$, and $V\subset U$ - a locally complete
intersection subvariety in $U$ of pure dimension $n-m$, locally defined in $U_{\alpha}$
by the holomorphic vector function
$${\bf F}^{(\alpha)}(z)=\left[\begin{tabular}{c}
$F^{(\alpha)}_1(z)$\vspace{0.1in}\\
\vdots\vspace{0.1in}\\
$F^{(\alpha)}_m(z)$
\end{tabular}\right],$$
i.e.
$$V\cap U_{\alpha}=\left\{z\in U_{\alpha}:\ F^{(\alpha)}_1(z)=\cdots=F^{(\alpha)}_m(z)=0\right\}.$$
%***************************************************************************************
%******* Definition  ConormalDualizing ***************************************************
\begin{definition}\label{ConormalDualizing} (Conormal and Dualizing Bundles)
The conormal vector bundle $N(V)$ on a locally complete intersection subvariety $V$ is defined by the
nondegenerate holomorphic transition matrices $A_{\alpha\beta}(z)\in H\left(U_{\alpha\beta}\right)$ such that
\begin{equation}\label{TransitionRelation}
{\bf F}^{(\alpha)}(z)=A_{\alpha\beta}(z)\cdot {\bf F}^{(\beta)}(z)
\end{equation}
on $U_{\alpha\beta}=U_{\alpha}\cap U_{\beta}$.\\
\indent
Following \cite{Gro} and \cite{Ha} we define the dualizing bundle on a locally complete intersection subvariety $V$ as
\begin{equation}\label{Dualizing}
\omega^{\circ}_V=\omega_{\C P^n}\otimes \det N(V)^{-1}
\end{equation}
where $\omega_{\C P^n}$ is the canonical bundle on $\C P^n$.
\end{definition}
%***************************************************************************************
%***************************************************************************************
{\bf Remark.} Adjunction formula (see Proposition 8.20 in Ch. II of \cite{Ha}) shows that for a nonsingular
$V$ the bundle defined in \eqref{Dualizing} coincides with the canonical bundle $\omega_V$, implicitly
used in Theorem~\ref{OldTheorem}, making $\omega^{\circ}_V$ a natural generalization of the canonical bundle
for locally complete intersection subvarieties of $\C P^n$.

\indent
We define further the spaces of residual currents and of residual $\bar\partial$-cohomologies on $V_D$, where
$V\subset \C P^n$ is a locally complete intersection subvariety, and $D$ a domain in $\C P^n$. In what follows we
denote by ${\cal E}$ the space of infinitely differentiable functions.

%**************************************************************************
%******* Definition  Residual ***************************************************
\begin{definition}\label{Residual} (Residual Currents) For a subvariety $V\subset \C P^n$ of the pure
dimension $n-m$ locally satisfying \eqref{V_cap_U} we say that a $(n,m+q)$ current
$\phi$ with support in $V$ is a residual current $\phi\in C^{(0,q)}\left(V_D,\omega_V^{\circ}\right)$
if there exists a finite collection of open neighborhoods $\left\{U_{\alpha}\subset \C P^n\right\}_{\alpha=1}^N$
and differential forms $\Phi_{\alpha}\in {\cal E}^{(n,q)}(U_{\alpha}\cap D)$, such that
\begin{equation}\label{ResCurrent}
\begin{cases}
\bigcup_{\alpha=1}^NU_{\alpha}\supset V,\\
{\dis \langle\phi,\psi\rangle=\int_{U_{\alpha}}\psi\wedge\Phi_{\alpha}
\wedge\bar\partial\left(\frac{1}{F^{(\alpha)}_1}\right)\wedge\cdots
\wedge\bar\partial\left(\frac{1}{F^{(\alpha)}_m}\right)
\stackrel{\rm def}{=}\lim_{t\to 0}\int_{T^{\epsilon}_{\left\{{\bf F}^{(\alpha)}\right\}}(t)}
\frac{\psi\wedge\Phi_{\alpha}}{\prod_{k=1}^m F^{(\alpha)}_k}, }\\
\Phi_{\alpha}=\left(\det{A_{\alpha\beta}}\right)^{-1}\cdot\Phi_{\beta}
+\sum_{k=1}^mF^{(\alpha)}_k\cdot\Omega^{(\alpha\beta)}_k\ \mbox{on}\ U_{\alpha}\cap U_{\beta}\cap D,
\end{cases}
\end{equation}
where $\psi\in {\cal E}_c^{(0,n-m-q)}(U_{\alpha}\cap D)$ is a smooth form with compact support in $U_{\alpha}\cap D$,
$$T^{\epsilon}_{\left\{{\bf F}^{(\alpha)}\right\}}(t)=\left\{|F^{(\alpha)}_1(z)|=\epsilon_1(t),
\dots,|F^{(\alpha)}_m(z)|=\epsilon_m(t)\right\}$$
is a family of tubular varieties depending on the real parameter $t$, the limit in the
right-hand side is taken along an admissible path $\left\{\epsilon_k(t)\right\}_1^m$ in
the sense of Coleff-Herrera-Lieberman \cite{CH}, \cite{HL}, i.e. an analytic map $\epsilon:[0,1]\to \R^m$ satisfying
the conditions
\begin{equation}\label{admissible}
\begin{cases}
\lim_{t\to 0}\epsilon_m(t)=0,\\
{\dis \lim_{t\to 0}\frac{\epsilon_j(t)}{\epsilon^l_{j+1}(t)}=0,\ \mbox{for any}\ l\in\Z },
\end{cases}
\end{equation}
$A_{\alpha\beta}$ are holomorphic matrices from \eqref{TransitionRelation},
and $\Omega^{(\alpha\beta)}_k\in {\cal E}\left(\ U_{\alpha}\cap U_{\beta}\cap D\right)$.\\
\indent
A residual current $\phi\in C^{(0,q)}\left(V_D,\omega_V^{\circ}\right)$ is called $\bar\partial$-closed -
$\phi\in Z^{(0,q)}\left(V_D,\omega_V^{\circ}\right)$, if the following condition is satisfied
\begin{equation}\label{Closedness}
\bar\partial\Phi_{\alpha}=\sum_{k=1}^mF^{(\alpha)}_k\cdot\Omega^{(\alpha)}_k\ \mbox{on}\ U_{\alpha}\cap D,
\end{equation}
where $\Omega^{(\alpha)}_k\in {\cal E}\left(\ U_{\alpha}\cap D\right)$.
\end{definition}
%************************************************************************************

{\bf Remarks.}
\begin{itemize}
\item
Condition \eqref{admissible}, though looking technical, can not be replaced by a simpler condition
$\epsilon_j(t)\to 0$, $t\to 0$, $j=1,\ldots,m$, as was shown by Passare and Tsikh in \cite{PT}.
\item
Notation in the definition above is substantiated by the fact that the collection
$$\left\{\Phi_{\alpha}\wedge\bar\partial\left(\frac{1}{F^{(\alpha)}_1}\right)\wedge\cdots
\wedge\bar\partial\left(\frac{1}{F^{(\alpha)}_m}\right)\right\}_{\alpha=1}^N$$
naturally defines a current of type $(0,q)$ on  $V_D$ with coefficients in holomorphic bundle
$\omega_V^{\circ}$ defined in \eqref{Dualizing}.
\end{itemize}

%************************************************************************************
%******* Definition  BarExact *************************************************************
\begin{definition}\label{BarExact} (Residual $\bar\partial$-cohomologies)
A $\bar\partial$-closed residual current
$\phi\in Z^{(0,q)}\left(V_D,\omega_V^{\circ}\right)$ is called $\bar\partial$-exact
($\phi\in B^{(0,q)}\left(V_D,\omega_V^{\circ}\right)$) if there exists a residual current
$\psi\in C^{(0,q-1)}(V_D,\omega_V^{\circ})$ such that $\bar\partial\psi=\phi$.\\
\indent
Therefore
$$B^{(0,q)}\left(V_D,\omega_V^{\circ}\right)\subseteq Z^{(0,q)}\left(V_D,\omega_V^{\circ}\right),$$
and the spaces of residual $\bar\partial$-cohomologies of $V_D$ of the type $(0,q)$:
$$H^{(0,q)}\left(V_D,\omega_V^{\circ}\right)=Z^{(0,q)}\left(V_D,\omega_V^{\circ}\right)/B^{(0,q)}\left(V_D,\omega_V^{\circ}\right)$$
are well defined.
\end{definition}
%*****************************************************************************************

\indent
Before defining the complex Radon transform we introduce an additional notation.
We denote by $S_V$ the following set of hyperplanes
$$S_V=\left\{\xi\in D^*:\ \dim_{\C}\left(V\cap\C P^{n-1}_{\xi}\right)\neq n-m-1\right\}.$$
Using the arguments similar to those in the proof of Bertini's theorem (see \cite{Ha}) we obtain that $S_V$
is a subset of an analytic set in $D^*$.

\indent
In the definitions below we define the complex Radon transform of residual currents and the Fantappi\'e transform
of linear functionals on $H^0\left(V,{\cal O}/{\cal I}\right)^{\prime}$.

%************************************************************************************
%******* Definition of Radon transform  ********************************************
\begin{definition}\label{RadonTransform} (Complex Radon Transform)  Let $V\subset \C P^n$ be a locally
complete intersection subvariety of pure dimension $n-m$. Then we define the Radon transform
$${\cal R}_V:Z^{(0,n-m-1)}\left(V_D,\omega_V^{\circ}\right)\to H^{(1,0)}\left(D^*\setminus S_V\right)$$
on the space of $\bar\partial$-closed residual currents by the formula (see Proposition~\ref{Coincidence})
\begin{multline}\label{Radon}
{\cal R}_V[\phi](\xi)=\frac{1}{\left(2\pi i\right)^{m+1}}
\sum_{j=0}^n\Bigg(\sum_{\alpha=1}^N\int_{D}
\vartheta_{\alpha}(z)\cdot z_j\cdot\Phi^{(n,n-m-1)}_{\alpha}(z)\\
\wedge\bar\partial\left(\frac{1}{\langle\xi\cdot z\rangle}\right)\bigwedge_{k=1}^m
\bar\partial\left(\frac{1}{F^{(\alpha)}_k(z)}\right)\Bigg)d\xi_j,
\end{multline}
where $\left\{\vartheta_{\alpha}\right\}_1^N$ is a partition of unity subordinate to a finite
cover $\left\{U_{\alpha}\right\}_1^N$ of $D$ by open subdomains in $\C P^n$,
and the forms
$$\left\{\Phi^{(n,n-m-1)}_{\alpha}\bigwedge_{k=1}^m
\bar\partial\left(\frac{1}{F^{(\alpha)}_k}\right)\right\}_{\alpha=1}^N$$
are the local representatives of the current $\phi$.
\end{definition}
%*****************************************************************************************

%************************************************************************************
%******* Definition of Fantappie transform  ********************************************
\begin{definition}\label{FantappieTransform} (Fantappi\'e Transform)  Let $V\subset \C P^n$ be a locally
complete intersection subvariety, let $G$ be a linearly convex compact in $\C P^n$, and let ${\cal I}$
be the sheaf of ideals, associated with $V$. We define the Fantappi\'e transform of a linear functional
$\mu\in H^0(G,{\cal O}/{\cal I})^{\prime}$ by the formula
\begin{equation}\label{Fantappie}
{\cal F}_V[\mu](\xi)=\sum_{j=0}^n\mu\left(\frac{z_j}{\langle \xi\cdot z\rangle}\right)d\xi_j,
\end{equation}
where $\xi\in D^*=({\C}P^n\setminus G)^*$.
\end{definition}
%*****************************************************************************************

\indent
The theorem below is the main result of the present article. In this theorem we describe the action
of the Fantappi\'e and complex Radon transforms on the spaces of residual cohomologies of linearly concave locally
complete intersection subvarieties of $\C P^n$.

%********************************************************************************
%******* Theorem RadonAction  ***************************************************
\begin{theorem}\label{RadonAction} Let
\begin{equation}\label{V-Equations}
V=\left\{z\in \C P^n:\ P_1(z)=\cdots=P_r(z)=0\right\}
\end{equation}
be a locally complete intersection subvariety of pure dimension $(n-m)$ with the structure sheaf ${\cal O}/{\cal I}$,
where ${\cal I}$ is the sheaf of ideals defined by homogeneous polynomials $\left\{P_k\right\}_1^r$, $r\geq m$.
Let $D\subset \C P^n$ be a linearly concave domain, and let $D^*$ be  its dual domain.\\
\indent
Then transform ${\cal R}_V$ defined in \eqref{Radon} induces a continuous
linear operator on the space of residual $\bar\partial$-cohomologies
$${\cal R}_V:H^{(0,n-m-1)}\left(V_D,\omega_V^{\circ}\right)\to H^{(1,0)}\left(D^*\right),$$
and transform ${\cal F}_V$ defined in \eqref{Fantappie} induces a continuous
linear operator
$${\cal F}_V:\ H^0(G,{\cal O}/{\cal I})^{\prime}\to H^{(1,0)}(D^*).$$
\indent
Transforms ${\cal R}_V$ and ${\cal F}_V$ satisfy the following properties:
\begin{itemize}
\item[(i)] $\text{Ker}\,{\cal F}_V=\{0\}$,
$\text{Ker}\,{\cal R}_V\subset H^{(0,n-m-1)}\left(V_D,\omega_V^{\circ}\right)$
is finite-dimensional  and consists of restrictions to $V_D$ of classes of residual $\bar\partial$-cohomologies from
$H^{(0,n-m-1)}\left(V,\omega_V^{\circ}\right)$,

\item[(ii)]  the images of ${\cal F}_V$ and ${\cal R}_V$ are the following subspaces in $H^{(1,0)}(D^*)$:
\begin{multline}\label{Images}
\text{Image}\ {\cal F}_V\\
=\left\{f\in H^{(1,0)}(D^*): f=dg\ \mbox{with}\ g\in H^0(D^*)\ \mbox{such that}\
\left\{P_k\left(\frac{\partial}{\partial\xi}\right)g=0\right\}_1^r\right\},\\
\text{Image}\ {\cal R}_V=\left\{f\in \text{Image}\ {\cal F}_V:\ f={\cal F}_V[\mu],\ \mbox{where}\
\mu(h)=0\ \text{for}\ \forall  h\in H^0\left(\C P^n,{\cal O}/{\cal I}\right)\right\},
\end{multline}

\item[(iii)] if $V$ is connected in the sense that $\dim H^0(V,{\cal O}/{\cal I})=1$, then
$$\text{Image}\ {\cal R}_V=\text{Image}\ {\cal F}_V,$$

\item[(iv)] for a functional $\mu\in H^0(G,{\cal O}/{\cal I})^{\prime}$ defined for $h\in H^0(G,{\cal O}/{\cal I})$
through the residual  current $\phi=\{\Phi_{\alpha}\}$ by the formula
$$\mu(h)=\sum_{\alpha=1}^N\int\limits_{bD_{\delta}}\vartheta_{\alpha}(z)
h(z)\Phi_{\alpha}(z)\bigwedge_{k=1}^m
\bar\partial\left(\frac{1}{F^{(\alpha)}_k(z)}\right),$$
the following equality holds
\begin{equation}\label{R-VEquality}
{\cal R}_V[\varphi](\xi)=\left(\frac{1}{2\pi i}\right)^{m+1}{\cal F}_V[\mu](\xi).
\end{equation}
\end{itemize}
\end{theorem}
%******************End Theorem RadonAction*************************************************
{\bf Remarks.}
\begin{itemize}
\item
The statements in (ii) of Theorem~\ref{RadonAction} can be interpreted as versions of the Ehrenpreis
\lq\lq fundamental principle\rq\rq\ for systems of partial differential equations (see \cite{E}, \cite{P1}, \cite{G})
in terms of Fantappi\'e and complex Radon transforms instead of Fourrier-Laplace transform.
\item
If $V$ is an arbitrary, not necessarily reduced, complete intersection in $\C P^n$ and $D$ is a linearly concave
domain in $\C P^n$, then in \cite{HP2} an explicit inversion formula for Radon transform ${\cal R}_V$ is obtained
together with a formula for solutions of appropriate boundary value problem for the corresponding system of
homogeneous differential equations with constant coefficients in $D^*$.
\item
For the case $m=n-1$ the statement (i) of Theorem~\ref{RadonAction} for Radon transform follows
from the result of Fabre \cite{F}.
\item
If $V$ is a complete intersection in $\C P^n$, then the property of $V$ to be connected in the sense of (iii)
is always satisfied (see Ex. 5.5 \S III.5 in \cite{Ha}).
\item
Theorem~\ref{RadonAction} admits a generalization for analytic subvarieties of a linearly concave
domain $D$. If $m<n-1$, then an analytic subvariety $V^{\prime}\subset D$ of $D$ is a trace of an algebraic
subvariety $V\subset \C P^n$ (see \cite{R}, \cite{Siu}), and an appropriate version of
Theorem~\ref{RadonAction} applies.
If $m=n-1$, then $V^{\prime}\subset D$ is a trace of an algebraic subvariety $V\subset \C P^n$ if there exists a form
$\phi\in  Z^{(0,n-m-1)}\left(V^{\prime},\omega_V^{\circ}\right)$, such that
$\phi\neq 0$ almost everywhere on $V$ and $R_{V^{\prime}}[\phi]\equiv 0$ (see \cite{Gr}, \cite{F}).
\end{itemize}

\indent
In section \ref{Properties} we prove the correctness of definition \ref{RadonTransform} and some properties
of ${\cal R}_V$ and ${\cal F}_V$, and in sections \ref{KernelDescription} and \ref{Images} we prove propositions representing different
parts of Theorem~\ref{RadonAction}.

\section{Properties of residual currents.}
\label{Properties}

\indent
In this section we describe some properties of residual currents used in the proof of Theorem~\ref{RadonAction}
and prove some properties of the Radon transform defined by formula \eqref{Radon}.
In the proposition below we describe the dependence of a local formula for a residual current on the choice of a basis
of the ideal for the case of a complete intersection.

%******************************************************************************************
%*******  Proposition TransformationLaw ***************************************************
\begin{proposition}\label{TransformationLaw}
Let $U\in \C P^n$ be a domain in $\C P^n$ and let $V\subset U$ be a complete
intersection subvariety of pure dimension $n-m$ in $U$, defined by two different collections
of holomorphic functions ${\bf F}=\left\{F_k\right\}_1^m$ and
${\bf P}=\left\{P_k\right\}_1^m$ such that
\begin{equation}\label{F-A-P}
{\bf F}=A\cdot{\bf P}
\end{equation}
where $A(z)$ is a nondegenerate holomorphic matrix-function.\\
\indent
Let $\left\{\epsilon(t)\right\}$ be an admissible path, and let
$$\begin{array}{ll}
T^{\epsilon}_{\left\{{\bf F}\right\}}(t)
=\left\{z\in U:\ |F_1(z)|=\epsilon_1(t),\dots,|F_m(z)|=\epsilon_m(t)
\right\},\vspace{0.1in}\\
T^{\epsilon}_{\left\{{\bf P}\right\}}(t)
=\left\{z\in U:\ |P_1(z)|=\epsilon_1(t),\dots,|P_m(z)|=\epsilon_m(t)
\right\},\end{array}$$
be the corresponding tubular varieties.\\
\indent
Then for an arbitrary $\gamma\in {\cal E}_c^{(n,n-m)}(U)$ we have the following equality
\begin{equation}\label{TransformationEquality}
\lim_{t\to 0}\int_{T^{\epsilon}_{\left\{{\bf P}\right\}}(t)}
\frac{\gamma(z)}{\prod_{k=1}^m P_k(z)}
=\lim_{t\to 0}\int_{T^{\epsilon}_{\left\{{\bf F}\right\}}(t)}
\frac{\det{A}(z)\cdot\gamma(z)}{\prod_{k=1}^m F_k(z)}.
\end{equation}
\end{proposition}
%**************************************************************************************
%***************************************************************************************
\indent
{\bf Proof.}\ In the proof of Proposition~\ref{TransformationLaw} we will use the
following proposition describing the transformation of the Grothendieck's residue
under the change of basis in the ideal for the case of isolated point in $\C^n$.

%******************************************************************************************
%*******  Proposition PointTransformationLaw *************************************************
\begin{proposition}\label{PointTransformationLaw}
(\cite{T},\cite{GrH})\ Let $U\in \C^n$ be a neighborhood
of the origin $\{0\}\in\C^n$ and let ${\bf P}= \{P_1,\dots,P_n\}$
and ${\bf F}= \{F_1,\dots,F_n\}$ be two different
collections of holomorphic functions on $U$ having $\{0\}$
as an isolated zero, and satisfying \eqref{F-A-P}
with a nondegenerate holomorphic matrix-function $A(z)$ on $U$.\\
\indent
Then for an arbitrary function $h\in {\cal E}_c(U)$ we have the following equality
\begin{equation}\label{DeterminantTransformation}
\lim_{t\to 0}\int_{T^{\epsilon}_{\left\{{\bf P}\right\}}(t)}
\frac{h(z)}{\prod_{k=1}^n P_k(z)}
=\lim_{t\to 0}\int_{T^{\epsilon}_{\left\{{\bf F}\right\}}(t)}
\frac{\det{A}(z)\cdot h(z)}{\prod_{k=1}^n F_k(z)}.
\end{equation}
\end{proposition}
\qed
%*******  End Proposition PointTransformationLaw ***************************************

\indent
To prove equality \eqref{TransformationEquality} we use the fibered residual currents
from \cite{CH}. Namely, we consider a polydisk ${\cal P}^n=\left\{|z_i|<1,\ i=1,\dots,n\right\}
\subset U$ such that the restriction of the projection
$$\pi:\ {\cal P}^n\to {\cal P}^{n-m},$$
defined by the formula $\pi(z_1,\dots,z_n)=(z_{m+1},\dots,z_n)$,
to $V\cap{\cal P}$ is a finite proper covering. Then we use Theorem 1.8.3 from \cite{CH}
and obtain the existence of a holomorphic function $g$ on ${\cal P}^n$ such that
$\dim\big\{V\cap\{|g(z)|=0\}\big\}\leq n-m-1$,
and
\begin{equation}\label{FiberedEquality}
\lim_{t\to 0}\int_{T^{\epsilon}_{\left\{{\bf P}\right\}}(t)}
\frac{\gamma(z)}{\prod_{k=1}^m P_k(z)}
=\lim_{\delta\to 0}\int_{V\cap\{|g(z)|>\delta\}}\mbox{res}_{\{{\bf P},\pi\}}
\left(\gamma,z\right),
\end{equation}
where
$$\mbox{res}_{\{{\bf P},\pi\}}\left(\gamma,z\right)
=\lim_{t\to 0}\int_{T^{\epsilon}_{\left\{{\bf{\widetilde P}}\right\}}(t)}
\frac{{\widetilde\gamma}\left(z_{m+1},\dots,z_n\right)}{\prod_{k=1}^n {\widetilde P}_k(z)},$$
$${\widetilde\gamma}\left(z_{m+1},\dots,z_n\right)=\gamma\Big|_{\pi^{-1}(z_{m+1},\dots,z_n)},
\hspace{0.05in}{\widetilde P}_k=P_k\Big|_{\pi^{-1}(z_{m+1},\dots,z_n)},$$
and $z\in V\cap \pi^{-1}(z_{m+1},\dots,z_n)$.\\
\indent
Applying Proposition~\ref{PointTransformationLaw}
to the right-hand side of equality \eqref{FiberedEquality} we obtain equality
\eqref{TransformationEquality}.
\qed
%*******  End Proposition TransformationLaw ***********************************************

\indent
The next proposition is a reformulation of Theorem 1.7.6(2) from \cite{CH}, which will be used in the article.
%********************************************************************************
%******* Proposition 1.7.6(2)  ***************************************************
\begin{proposition}\label{1.7.6(2)}
Let $U$ be a domain in $\C^n$, and let
$$V=\left\{z\in U:\ F_1(z)=\cdots=F_m(z)=0\right\}$$
be a complete intersection in $U$. If a differential form $\Phi\in {\cal E}_c^{(n,n-m)}\left(U\right)$
with compact support in $U$ admits a representation
$$\Phi=\sum_{k=1}^m F_k\cdot \Phi_k,$$
where forms $\Phi_k\in {\cal E}^{(n,n-m)}\left(U\right)$ have compact support in $U$, then
$$\int_U\Phi\bigwedge_{k=1}^m\bar\partial\left(\frac{1}{F_k}\right)=0.$$
\qed
\end{proposition}
%********************************************************************************
%*******  End Proposition 1.7.6(2)  ***********************************************

\indent
In the proposition below we prove the existence of a family of smoothly bounded linearly concave domains
approximating $D$. Existence of such family provides a convenient tool in many constructions of the present article.

%***********************************************************************************
%******* Proposition Exhaustion  *********************************************************
\begin{proposition}\label{Exhaustion}
Let a linearly concave domain $D\subset \C P^n$ admit a continuos map $\eta: D\to D^*$ satisfying
condition $\langle \eta(z)\cdot z\rangle=0$. Then there exist a sequence of real numbers
$\left\{\delta_n\right\}_1^{\infty}$ such that $\delta_n>\delta_m$ for $n<m$ and $\lim_{n\to\infty}\delta_n=0$,
and of smoothly bounded linearly concave domains
\begin{equation}\label{DeltaDomain}
D_{\delta_n}\subset D=\left\{z\in D: \rho_{\delta_n}(z)<0\right\}
\end{equation}
satisfying
\begin{equation}\label{DeltaDomains}
D_{\delta_n}\subset D_{\delta_m}\quad \text{for}\ m>n,\quad
\text{and}\quad \bigcup_{n=1}^{\infty} D_{\delta_n}=D.
\end{equation}
\end{proposition}
%**********************************************************************************
%**********************************************************************************
\indent
{\bf Proof.} We construct a sequence of smoothly bounded linearly concave domains satisfying \eqref{DeltaDomains}
in two steps. On the first step we construct a family of domains exhausting $D^*$.
We consider the function $\rho^*(\xi)={\rm dist}(\xi, bD^*)$ on $D^*$, and averaging this function with the kernel
$K_{\delta}(\zeta)=\delta^{-2n}\cdot K(\zeta/\delta)$, where
$$K(\zeta)=\begin{cases}
\begin{aligned}
&Ce^{1/(|\zeta|^2-1)}\ &\text{if}\ |\zeta|<1,\\
&0\ &\text{if}\ |\zeta|\geq 1,
\end{aligned}
\end{cases}$$
and $C=\left(\int_{|\zeta|\leq 1}e^{1/(|\zeta|^2-1)}d\zeta\right)^{-1}$, obtain a smooth function
\begin{equation*}
\rho^*_{\delta}(\xi)=\int\rho^*(\zeta)K_{\delta}(\xi-\zeta)d\zeta
\end{equation*}
on the set $\left\{\xi\in D^*:\rho^*(\xi)>\delta\right\}$. We define then for $\nu<\delta/2$
$$D^*_{\delta,\nu}=\left\{\xi\in D^*: \rho^*_{\delta}(\xi)>3\delta-\nu\right\}.$$
\indent
To see that
\begin{equation}\label{Embedding}
\left\{\xi\in D^*:\rho^*(\xi)>4\delta\right\}\subset
D^*_{\delta,\nu}\subset \left\{\xi\in D^*:\rho^*(\xi)>\delta\right\}
\end{equation}
for $\nu<\delta/2$ we use the inequality
\begin{multline*}
\left|\rho^*_{\delta}(\xi)-\rho^*(\xi)\right|=\left|\int\left(\rho^*(\zeta)-\rho^*(\xi)\right)
K_{\delta}(\xi-\zeta)d\zeta\right|\\
\leq\delta^{-2n}\int\left|\rho^*(\zeta)-\rho^*(\xi)\right|K\left(\frac{\xi-\zeta}{\delta}\right)d\zeta
=\int_{|u|\leq 1}\left|\rho^*(\xi+\delta\cdot u)-\rho^*(\xi)\right|K(u)du\leq\delta.
\end{multline*}
Relation \eqref{Embedding} shows that the family of domains $D^*_{\delta,\nu}$ exhausts domain $D^*$.\\
\indent
On the second step we consider the domain
$$W^*_{\delta}=\left\{\xi\in D^*:\rho^*(\xi)>\delta\right\},$$
and apply a smoothing procedure, similar to the described above, to the  continuous family of hyperplanes
$\eta:D\to D^*$ restricted to the domain $\eta^{-1}\left(W^*_{\delta}\right)$.
For $z$ in the domain
$$\eta^{-1}\left(W^*_{\delta}\right)\cap U_j=\left\{z\in \eta^{-1}\left(W^*_{\delta}\right):
z_j\neq 0\right\}$$
we define
$$\eta^i_{j,\delta^{\prime}}(z)=\begin{cases}
{\dis \int\eta_i(\zeta)K_{\delta^{\prime}}(z-\zeta)d\zeta\quad \text{if}\ i\neq j, }\vspace{0.1in}\\
{\dis \eta^j_{j,\delta^{\prime}}(z)=-\sum_{i\neq j}\eta^i_{j,\delta^{\prime}}(z)\frac{z_i}{z_j}, }
\end{cases}$$
for $\delta^{\prime}>0$ small enough, and set
$$\eta_{\delta^{\prime}}(z)=\left(\eta_{0,\delta^{\prime}}(z),\dots,\eta_{n,\delta^{\prime}}(z)\right),$$
where
$$\eta_{k,\delta^{\prime}}(z)=\sum_{j=0}^n\vartheta_j(z)\cdot \eta^k_{j,\delta^{\prime}}(z),$$
and $\left\{\vartheta_j\right\}_{j=0}^n$ is a partition of unity subordinate
to the cover $\left\{U_j\right\}$ of $\C P^n$.\\
\indent
We notice that
for every $j\in (0,\dots,n)$ we have
$$\sum_{k=0}^nz_k\cdot\eta^k_{j,\delta^{\prime}}(z)=0,$$
and therefore
$$\sum_{k=0}^nz_k\cdot\eta_{k,\delta^{\prime}}(z)=\sum_{k=0}^nz_k\cdot
\left(\sum_{j=0}^n\vartheta_j(z)\cdot\eta^k_{j,\delta^{\prime}}(z)\right)$$
$$=\sum_{j=0}^n\vartheta_j(z)\cdot\left(\sum_{k=0}^nz_k\cdot\eta^k_{j,\delta^{\prime}}(z)\right)=0.$$
\indent
Then we obtain a continuous and smooth in a neighborhood of $\eta^{-1}\left(D^*_{\delta}\right)$
family of hyperplanes $\eta_{\delta^{\prime}}(z)\in D^*$ such that $z\in\eta_{\delta^{\prime}}(z)$ for
$$z\in \eta^{-1}\left\{\xi\in D^*: \rho^*(\xi)>\delta\right\}.$$
\indent
We define
$$D^{\prime}_{\delta,\nu}=\left\{z\in D: \C P^{n-1}(z)\subset D^*_{\delta,\nu}\right\}
=\left\{z\in D:\rho_{\delta}(z)\stackrel{\rm def}{=}
3\delta-\nu-\rho^*_{\delta}\big(\eta_{\delta^{\prime}}(z)\big)<0\right\},$$
and applying the Sard's theorem find $\nu^{\prime}<\delta/2$ such that
$D_{\delta}=D^{\prime}_{\delta,\nu^{\prime}}$ has smooth boundary.\\
\indent
Sequences $\left\{\delta_n\right\}_1^{\infty}$ and $\left\{D_{\delta_n}\right\}_1^{\infty}$
satisfying \eqref{DeltaDomains} can be chosen as subsequences corresponding to an arbitrary sequence
of decreasing $\delta_n$ tending to zero as $n\to\infty$ based on the exhaustion property. To construct an
\lq\lq explicit\rq\rq\ sequence  $\left\{D_{\delta_n}\right\}_1^{\infty}$ satisfying \eqref{DeltaDomains}
we can choose for example the sequence $\left\{\delta_n=\delta_1/8^{n-1}\right\}_1^{\infty}$.
The numbers $\delta^{\prime}_n$ can be chosen so that
\begin{equation}\label{rhoEstimate}
\left|\rho^*_{\delta_n}(\eta_{\delta^{\prime}_n}(z))-\rho^*_{\delta_n}(\eta(z))\right|<\frac{\delta_n}{16}
\end{equation}
and therefore $\left|\rho^*_{\delta_n}(\eta_{\delta^{\prime}}(z))-\rho^*(\eta(z))\right|\leq\frac{17}{16}\delta_n$,
for $z\in \eta^{-1}\left\{\xi\in D^*: 4\delta_n>\rho^*(\xi)>\delta_n\right\}$.\\
\indent
The boundary of the domain
$$D^{\prime}_{\delta_n,\nu_n}=\left\{z\in D:3\delta_n-\nu_n
-\rho^*_{\delta_n}\big(\eta_{\delta^{\prime}_n}(z)\big)<0\right\}$$
will satisfy the condition $\rho^*_{\delta_n}\big(\eta_{\delta^{\prime}_n}(z)\big)=3\delta_n-\nu_n$,
and for $z\in bD^{\prime}_{\delta_n,\nu_n}$ we will have using \eqref{rhoEstimate}
$$\frac{57}{16}\delta_n\geq\rho^*(\eta(z))\geq \frac{21}{16}\delta_n.$$
\indent
Since $\frac{57}{16\cdot 8}<\frac{21}{16}$ we obtain that
$bD^{\prime}_{\delta_n,\nu_n}\cap bD^{\prime}_{\delta_{n+1},\nu_{n+1}}=\emptyset$, and therefore
the sequence $D_{\delta_n}$ is strictly monotonous.
\qed\\
%******* End Proposition Exhaustion  *********************************************************

\indent
In the next proposition we prove a useful boundary formula for the Radon transform.
As a corollary of this formula we obtain that definition \eqref{Radon} of the Radon transform ${\cal R}_V$ coincides
with the standard definition of Radon transform for the case of a differential form on a nonsingular variety $V$.

%********************************************************************************
%******* Proposition Coincidence  ***************************************************
\begin{proposition}\label{Coincidence}
Let $D\subset \C P^n$ be a linearly concave domain,
let $\left\{U_{\alpha}\right\}_{\alpha=1}^N$ be a finite cover of $D$, let
$\left\{\vartheta_{\alpha}\right\}_{\alpha=1}^N$  be a partition of unity subordinate to the cover $\left\{U_{\alpha}\right\}_{\alpha=1}^N$, and
let $V\subset D$ be a locally complete intersection subvariety
of pure dimension $n-m$, locally defined in $U_{\alpha}$ by the holomorphic functions
$\left\{F^{(\alpha)}_k\right\}_{k=1}^m$.\\
\indent
Then for a $\bar\partial$-closed residual current $\phi$ defined locally by the differential forms
$$\Phi_{\alpha}\in {\cal E}^{(n,n-m-1)}\left(U_{\alpha}\right)$$
and a subdomain $D_{\delta}\subset D$ with smooth boundary $bD_{\delta}$ the following equality holds
\begin{multline}\label{RadonFantappie}
\sum_{\alpha=1}^N\int_{D}
\vartheta_{\alpha}(z)\cdot z_j\cdot\Phi_{\alpha}(z)
\wedge\bar\partial\left(\frac{1}{\langle\xi\cdot z\rangle}\right)\bigwedge_{k=1}^m
\bar\partial\left(\frac{1}{F^{(\alpha)}_k(z)}\right)\\
=\sum_{\alpha=1}^N\lim_{\tau\to 0}\lim_{t\to 0}
\int_{T^{\epsilon}_{\left\{{\bf F}^{(\alpha)},\tau\right\}}(t)}
\vartheta_{\alpha}(z)\frac{z_j\cdot\Phi_{\alpha}(z)}
{\langle\xi\cdot z\rangle\cdot\prod_{k=1}^m F^{(\alpha)}_k(z)}\\
=\sum_{\alpha=1}^N\lim_{t\to 0}
\int_{bD_{\delta}\cap T^{\epsilon}_{\left\{{\bf F}^{(\alpha)}\right\}}(t)}
\vartheta_{\alpha}(z)\frac{z_j\cdot\Phi_{\alpha}(z)}
{\langle\xi\cdot z\rangle\cdot\prod_{k=1}^m F^{(\alpha)}_k(z)},
\end{multline}
where
$$T^{\epsilon}_{\left\{{\bf F}^{(\alpha)},\tau\right\}}(t)
=\left\{z\in U_{\alpha}:\ \left\{|F^{(\alpha)}_k(z)|=\epsilon_k(t)\right\}_{k=1}^m,\ 
\chi(\xi,z)\stackrel{\rm def}{=}
\sum_{\alpha=1}^N\vartheta_{\alpha}(z)\cdot\left|\langle\xi\cdot z^{(\alpha)}\rangle\right|
=\tau\right\}$$
with admissible path $\left\{\epsilon_k(t)\right\}_{k=1}^m$.\\
\indent
Under the hypotheses of Theorem~\ref{RadonAction} the transform ${\cal R}_V$ from
\eqref{Radon} maps the $\bar\partial$-closed residual currents on $D$ with support on $V_D$
into holomorphic forms on $D^*$, and induces a linear map on the spaces of cohomologies.
\end{proposition}
%********************************************************************************
In the proof of Proposition~\ref{Coincidence} we will use the following two lemmas.

%***********************************************************************************
%*******  Lemma  ZeroIntegral **********************************************************
\begin{lemma}\label{ZeroIntegral} Let $U\subset \C P^n$ be a domain in $\C P^n$,
let $\left\{U_{\alpha}\right\}_{\alpha=1}^N$ be a finite cover of $U$, and let
$V\subset U$ be a locally complete intersection subvariety in $U$ of pure dimension $n-m$,
locally defined in $U_{\alpha}$ by holomorphic functions $\left\{F^{(\alpha)}_k\right\}_1^m$.
Let $\omega$ be a $\bar\partial$-closed residual current with support on $V$
locally defined by the differential forms $\Omega_{\alpha}\in {\cal E}^{(n,n-m-1)}(U_{\alpha})$.\\
\indent
Then for an arbitrary function $\eta\in {\cal E}_c\left(U\right)$ we have
\begin{equation}\label{ZeroEquality}
\sum_{\alpha=1}^N\lim_{t\to 0}\int_{T^{\epsilon}_{\left\{{\bf F}^{(\alpha)}\right\}}(t)}
\eta(z)\bar\partial\vartheta_{\alpha}(z)\wedge\frac{\Omega_{\alpha}(z)}
{\prod_{k=1}^m F^{(\alpha)}_k(z)}=0.
\end{equation}
\end{lemma}
%***************************************************************************************
%***************************************************************************************
\indent
{\bf Proof.} To prove equality \eqref{ZeroEquality} we apply the Stokes' formula,
and using equality
$$\bar\partial\Omega_{\alpha}=\sum_{k=1}^mF^{(\alpha)}_k\cdot\Omega^{(\alpha)}_k$$
for $i=1\dots N$ and Proposition~\ref{1.7.6(2)} obtain
\begin{equation*}
\sum_{\alpha=1}^N\lim_{t\to 0}\int_{T^{\epsilon}_{\left\{{\bf F}^{(\alpha)}\right\}}(t)}
\eta(z)\bar\partial\vartheta_{\alpha}(z)\wedge\frac{\Omega_{\alpha}(z)}
{\prod_{k=1}^m F^{(\alpha)}_k(z)}
\end{equation*}
\begin{multline*}
=-\sum_{\alpha=1}^N\lim_{t\to 0}\int_{T^{\epsilon}_{\left\{{\bf F}^{(\alpha)}\right\}}(t)}
\vartheta_{\alpha}(z)\bar\partial\eta(z)\wedge\frac{\Omega_{\alpha}(z)}
{\prod_{k=1}^m F^{(\alpha)}_k(z)}\\
-\sum_{\alpha=1}^N\lim_{t\to 0}\int_{T^{\epsilon}_{\left\{{\bf F}^{(\alpha)}\right\}}(t)}
\vartheta_{\alpha}(z)\eta(z)\wedge\frac{\bar\partial\Omega_{\alpha}(z)}
{\prod_{k=1}^m F^{(\alpha)}_k(z)}\\
=-\omega\left(\sum_{\alpha=1}^N\vartheta_{\alpha}\cdot\bar\partial\eta\right)
=-\omega\left(\bar\partial\eta\right)=\pm\bar\partial\omega(\eta)=0.
\end{multline*}
\qed
%*******  End Lemma ZeroIntegral ************************************************

%********************************************************************************
%******* Lemma TubeIndependence  ***************************************************
\begin{lemma}\label{TubeIndependence}
Let $D\subset \C P^n$ be a linearly concave domain,
let $V\subset D$ be a locally complete intersection subvariety
of pure dimension $n-m$, locally defined in $U_{\alpha}$ by holomorphic functions
$\left\{F^{(\alpha)}_k\right\}_1^m$.\\
\indent
Then for a fixed $\xi\in D^*\setminus S_V$ and a $\bar\partial$-closed residual current $\omega$
defined locally by the differential forms
$$\Omega_{\alpha}\in {\cal E}^{(n,n-m-1)}\left(U_{\alpha}\setminus
\left\{\langle\xi\cdot z\rangle=0\right\}\right)$$
the expression
\begin{equation}\label{TubeIntegral}
\sum_{\alpha=1}^N\lim_{t\to 0}
\int_{T^{\epsilon}_{\left\{{\bf F}^{(\alpha)},\tau\right\}}(t)}
\vartheta_{\alpha}(z)\frac{\Omega_{\alpha}(z)}{\prod_{k=1}^m F^{(\alpha)}_k(z)}
\end{equation}
is well defined, doesn't depend on $\tau$, and the following equality holds
\begin{equation}\label{BoundaryEquality}
\sum_{\alpha=1}^N\lim_{t\to 0}
\int_{T^{\epsilon}_{\left\{{\bf F}^{(\alpha)},\tau\right\}}(t)}
\vartheta_{\alpha}(z)\frac{\Omega_{\alpha}(z)}{\prod_{k=1}^m F^{(\alpha)}_k(z)}
=\sum_{\alpha=1}^N\lim_{t\to 0}
\int_{bD_{\delta}\cap T^{\epsilon}_{\left\{{\bf F}^{(\alpha)}\right\}}(t)}
\vartheta_{\alpha}(z)\frac{\Omega_{\alpha}(z)}{\prod_{k=1}^m F^{(\alpha)}_k(z)}.
\end{equation}
\end{lemma}
%********************************************************************************
\indent
{\bf Proof.} We fix a sufficiently small $\mu>0$ and consider for an arbitrary
$\tau>0$ such that $\tau<\mu$ a family of nonnegative functions
$\eta_{\nu}\in {\cal E}\left(D\right)$ such that
$$\eta_{\nu}(z)=\begin{cases}
0\ \mbox{if}\ \chi(\xi,z)<\tau-\nu,\\
1\ \mbox{if}\ \tau+\nu<\chi(\xi,z)<\mu,\\
0\ \mbox{if}\ \chi(\xi,z)>2\mu,
\end{cases}$$
and such that $\eta_{\nu}(z)=\eta(z)$ for $z$ with $\chi(\xi,z)>\mu$, where $1\geq\eta(z)\geq 0$
is a fixed smooth function.\\
\indent
Applying then the Stokes' formula we obtain
\begin{multline*}
\sum_{\alpha=1}^N\lim_{t\to 0}\int_{T^{\epsilon}_{\left\{{\bf F}^{(\alpha)}\right\}}(t)}
\eta_{\nu}(z)\bar\partial\vartheta_{\alpha}(z)\wedge\frac{\Omega_{\alpha}(z)}
{\prod_{k=1}^m F^{(\alpha)}_k(z)}\\
+\sum_{\alpha=1}^N\lim_{t\to 0}\int_{T^{\epsilon}_{\left\{{\bf F}^{(\alpha)}\right\}}(t)}
\vartheta_{\alpha}(z)\bar\partial\eta_{\nu}(z)\wedge\frac{\Omega_{\alpha}(z)}
{\prod_{k=1}^m F^{(\alpha)}_k(z)}\\
+\sum_{\alpha=1}^N\lim_{t\to 0}\int_{T^{\epsilon}_{\left\{{\bf F}^{(\alpha)}\right\}}(t)}
\vartheta_{\alpha}(z)\eta_{\nu}(z)\wedge\frac{\bar\partial\Omega_{\alpha}(z)}
{\prod_{k=1}^m F^{(\alpha)}_k(z)}=0,
\end{multline*}
which we transform using Proposition~\ref{1.7.6(2)} and  Lemma~\ref{ZeroIntegral} into
\begin{equation*}\label{StokesPhi}
\sum_{\alpha=1}^N\lim_{t\to 0}\int_{T^{\epsilon}_{\left\{{\bf F}^{(\alpha)}\right\}}(t)}
\vartheta_{\alpha}(z)\bar\partial\eta_{\nu}(z)\wedge\frac{\Omega_{\alpha}(z)}
{\prod_{k=1}^m F^{(\alpha)}_k(z)}=0,
\end{equation*}
and then further into
\begin{multline*}
-\sum_{\alpha=1}^N\lim_{t\to 0}\int_{T^{\epsilon}_{\left\{{\bf F}^{(\alpha)}\right\}}(t)
\cap\left\{\tau-\epsilon<\chi(\xi,z)<\tau+\epsilon\right\}}
\vartheta_{\alpha}(z)\bar\partial\eta_{\nu}(z)\wedge\frac{\Omega_{\alpha}(z)}
{\prod_{k=1}^m F^{(\alpha)}_k(z)}\\
=\sum_{\alpha=1}^N\lim_{t\to 0}\int_{T^{\epsilon}_{\left\{{\bf F}^{(\alpha)}\right\}}(t)
\cap\left\{\mu<\chi(\xi,z)<2\mu\right\}}
\vartheta_{\alpha}(z)\bar\partial\eta(z)\wedge\frac{\Omega_{\alpha}(z)}
{\prod_{k=1}^m F^{(\alpha)}_k(z)}
\end{multline*}
for arbitrary small $\tau>0$.\\
\indent
Considering then the limit of the equality above as $\nu\to 0$ we obtain the equality
\begin{multline}\label{EtaEquality}
\sum_{\alpha=1}^N\lim_{t\to 0}
\int_{T^{\epsilon}_{\left\{{\bf F}^{(\alpha)},\tau\right\}}(t)}
\vartheta_{\alpha}(z)\frac{\Omega_{\alpha}(z)}{\prod_{k=1}^m F^{(\alpha)}_k(z)}\\
=\sum_{\alpha=1}^N\lim_{t\to 0}\int_{T^{\epsilon}_{\left\{{\bf F}^{(\alpha)}\right\}}(t)
\cap\left\{\mu<\chi(\xi,z)<2\mu\right\}}
\vartheta_{\alpha}(z)\bar\partial\eta(z)\wedge\frac{\Omega_{\alpha}(z)}
{\prod_{k=1}^m F^{(\alpha)}_k(z)}.
\end{multline}
The limit in the right-hand side of \eqref{EtaEquality} exists according to the following proposition,
which is a reformulation of item (2) of Theorem 1.7.2 from \cite{CH}.

%********************************************************************************
%******* Proposition 1.7.2(2)  *********************************************************
\begin{proposition}\label{1.7.2(2)}
Let $U$ be a relatively compact domain in $\C^n$, let
$$V=\left\{z\in U: F_1(z)=\cdots=F_{m}(z)=0\right\}$$
be a complete intersection subvariety in $U$, let $\beta\in{\cal E}_c^{(n,n-m)}(U)$ be a differential
form with compact support in $U$, and let $T^{\epsilon}_{\left\{\bf F\right\}}(t)$ be an admissible path. Then
the following limit
$$\lim_{t\to 0}\int_{T^{\epsilon}_{\left\{{\bf F}\right\}}(t)}\frac{\beta(z)}{\prod_{k=1}^m F_k(z)}$$
exists.\qed
\end{proposition}
%********************************************************************************

From the form of the integral in the right-hand side of \eqref{EtaEquality} we conclude that it doesn't depend
on the choice of $\tau$, and therefore the same is true for the left-hand side of this equality.\\
\indent
To prove equality \eqref{BoundaryEquality} we change the choice of the family
of functions $\eta_{\nu}$ to the following:
$$\eta_{\nu}(z)=\begin{cases}
0\ \mbox{if}\ \chi(\xi,z)<\tau-\nu\ \mbox{or}\ \rho_{\delta}(z)>\nu,\\
1\ \mbox{if}\ \chi(\xi,z)>\tau+\nu\ \mbox{and}\ \rho_{\delta}(z)<-\nu,
\end{cases}$$
where $\rho_{\delta}$ is the function from Proposition~\ref{Exhaustion}.\\
\indent
Applying then the Stokes' formula we obtain the equality
\begin{multline*}
\sum_{\alpha=1}^N\lim_{t\to 0}\int_{T^{\epsilon}_{\left\{{\bf F}^{(\alpha)}\right\}}(t)}
\eta_{\nu}(z)\bar\partial\vartheta_{\alpha}(z)\wedge\frac{\Omega_{\alpha}(z)}
{\prod_{k=1}^m F^{(\alpha)}_k(z)}\\
+\sum_{\alpha=1}^N\lim_{t\to 0}\int_{T^{\epsilon}_{\left\{{\bf F}^{(\alpha)}\right\}}(t)}
\vartheta_{\alpha}(z)\bar\partial\eta_{\nu}(z)\wedge\frac{\Omega_{\alpha}(z)}
{\prod_{k=1}^m F^{(\alpha)}_k(z)}\\
+\sum_{\alpha=1}^N\lim_{t\to 0}\int_{T^{\epsilon}_{\left\{{\bf F}^{(\alpha)}\right\}}(t)}
\vartheta_{\alpha}(z)\eta_{\nu}(z)\wedge\frac{\bar\partial\Omega_{\alpha}(z)}
{\prod_{k=1}^m F^{(\alpha)}_k(z)}=0,
\end{multline*}
and then using Lemma~\ref{ZeroIntegral} and Proposition~\ref{1.7.6(2)}
\begin{multline*}
-\sum_{\alpha=1}^N\lim_{t\to 0}\int_{T^{\epsilon}_{\left\{{\bf F}^{(\alpha)}\right\}}(t)
\cap\left\{\tau-\nu<\chi(\xi,z)<\tau+\nu\right\}}
\vartheta_{\alpha}(z)\bar\partial\eta_{\nu}(z)\wedge\frac{\Omega_{\alpha}(z)}
{\prod_{k=1}^m F^{(\alpha)}_k(z)}\\
=\sum_{\alpha=1}^N\lim_{t\to 0}\int_{T^{\epsilon}_{\left\{{\bf F}^{(\alpha)}\right\}}(t)
\cap\left\{-\nu<\rho_{\delta}(z)<\nu\right\}}
\vartheta_{\alpha}(z)\bar\partial\eta_{\nu}(z)\wedge\frac{\Omega_{\alpha}(z)}
{\prod_{k=1}^m F^{(\alpha)}_k(z)}.
\end{multline*}
\indent
Passing to the limit as $\nu\to 0$ we obtain equality \eqref{BoundaryEquality}.
\qed\\
%*******  End Lemma TubeIndependence ******************************************

\indent
{\bf Proof of Proposition~\ref{Coincidence}.} Equality \eqref{RadonFantappie} is an immediate corollary of equality
\eqref{BoundaryEquality}. In view of \eqref{RadonFantappie} formula \eqref{Radon} defines a bounded
holomorphic function on the intersection of an arbitrary compact set in $D^*$ with $D^*\setminus S_V$.
Therefore, since $S_V$ is a subset of an analytic set in $D^*$,  there exists a unique extension of this function to $D^*$.\\
\indent
To prove that ${\cal R}_V[\phi]=0$ for a $\bar\partial$-exact residual current $\phi^{n,n-1}$
we assume the existence of a current
$\psi\in {\cal K}^{n,n-2}(D)$ such that equality
$$\left\langle\gamma^{(0,1)},\ \phi^{(n,n-1)}\right\rangle
=\left\langle\bar\partial\gamma^{(0,1)},\ \psi^{(n,n-2)}\right\rangle$$
is satisfied for an arbitrary $\gamma^{(0,1)}\in {\cal E}_c^{(0,1)}(D)$.\\
\indent
Then, using formula \eqref{RadonFantappie} and Proposition~\ref{1.7.6(2)} we obtain
\begin{multline*}
{\cal R}_V[\phi](\xi)=\frac{1}{(2\pi i)^{m+1}}
\sum_{j=0}^n\left(\sum_{\alpha=1}^N\lim_{t\to 0}
\int_{bD_{\delta}\cap T^{\epsilon}_{\left\{{\bf F}^{(\alpha)}\right\}}(t)}
\vartheta_{\alpha}(z)\frac{z_j\cdot\Phi_{\alpha}(z)}
{\langle\xi\cdot z\rangle\cdot\prod_{k=1}^m F^{(\alpha)}_k(z)}\right)d\xi_j\\
=\frac{1}{(2\pi i)^{m+1}}\sum_{j=0}^n\left(\sum_{\alpha=1}^N\lim_{\nu\to 0}\lim_{t\to 0}
\int_{T^{\epsilon}_{\left\{{\bf F}^{(\alpha)}\right\}}(t)}\vartheta_{\alpha}(z)
\bar\partial\eta_{\nu}(z)\wedge\frac{z_j\cdot\Phi_{\alpha}(z)}
{\langle\xi\cdot z\rangle\cdot\prod_{k=1}^m F^{(\alpha)}_k(z)}\right)d\xi_j\\
\end{multline*}
\begin{multline*}
+\frac{1}{(2\pi i)^{m+1}}\sum_{j=0}^n\left(\sum_{\alpha=1}^N\lim_{\nu\to 0}\lim_{t\to 0}
\int_{T^{\epsilon}_{\left\{{\bf F}^{(\alpha)}\right\}}(t)}\vartheta_{\alpha}(z)
\eta_{\nu}(z)\frac{z_j\cdot\bar\partial\Phi_{\alpha}(z)}
{\langle\xi\cdot z\rangle\cdot\prod_{k=1}^m F^{(\alpha)}_k(z)}\right)d\xi_j\\
=\frac{1}{(2\pi i)^{m+1}}\sum_{j=0}^n\left(\sum_{\alpha=1}^N\lim_{\nu\to 0}
\left\langle\frac{z_j}{\langle\xi\cdot z\rangle}\cdot\bar\partial\eta_{\nu},\
\phi\right\rangle\right)d\xi_j\\
=\frac{1}{(2\pi i)^{m+1}}\sum_{j=0}^n\left(\sum_{\alpha=1}^N\lim_{\nu\to 0}
\left\langle\frac{z_j}{\langle\xi\cdot z\rangle}\cdot\bar\partial^2\eta_{\nu},\
\psi\right\rangle\right)d\xi_j=0,
\end{multline*}
where
$$\eta_{\nu}(z)=\begin{cases}
1\ \mbox{if}\ \rho_{\delta}(z)>\nu,\\
0\ \mbox{if}\ \rho_{\delta}(z)<-\nu.
\end{cases}$$
\indent
The equality above allows to define the Radon transform ${\cal R}_V[\phi]=0$ for an arbitrary
$\bar\partial$-exact current $\phi$ with support on $V\cap D$.
\qed\\
%******* End Proposition Coincidence  ***************************************************

\indent
In the proposition below we prove the inclusion of the images of ${\cal F}_V$ and ${\cal R}_V$
in the space of solutions in the right-hand side of \eqref{Images}.

%***********************************************************************************
%******* Proposition F-RInclusion  **************************************************
\begin{proposition}\label{F-RInclusion} Radon and Fantappi\'e transforms defined in \eqref{Radon}
and \eqref{Fantappie} satisfy the following properties:
\begin{multline}\label{Inclusions}
\text{Image}\ {\cal F}_V\\
\subseteq \left\{f\in H^{(1,0)}(D^*): f=dg\ \mbox{with}\ g\in H^0(D^*)\ \mbox{such that}\
\left\{P_k\left(\frac{\partial}{\partial\xi}\right)g=0\right\}_1^r\right\},\\
\text{Image}\ {\cal R}_V\subseteq\left\{f\in \text{Image}\ {\cal F}_V:\ f={\cal F}_V[\mu],\ \mbox{where}\
\mu(h)=0\ \text{for}\ \forall  h\in H^0\left(\C P^n,{\cal O}/{\cal I}\right)\right\}.
\end{multline}
\end{proposition}
%***********************************************************************************
\indent
{\bf Proof.} For the Fantappi\'e transform of a linear functional $\mu\in H^0(G,{\cal O}/{\cal I})^{\prime}$
we have
$${\cal F}[\mu]=\sum_{j=0}^n\mu\left(\frac{z_j}{\langle \xi\cdot z\rangle}\right)d\xi_j
=d_{\xi}g(\xi),$$
where
\begin{equation}\label{gFantappie}
g(\xi)=\mu\left(\log\langle \xi\cdot z\rangle\right).
\end{equation}
\indent
We notice that analytic function $\log\left\langle\xi\cdot z\right\rangle$, and therefore $g(\xi)$, is well defined
on $D^*(z)$. It is a corollary of the contractibility of
$$D^*(z)=\left\{\xi\in D^*:  \langle\xi\cdot z\rangle=0\right\}$$
for any $z\in D$ under the condition of existence of a continuous family of hyperplanes covering the whole $D$.
Namely, as it was proved in \cite{GH}, the existence of such family implies the isomorphism
$$H(\C P^n\setminus D)^{\prime}\cong H(D^*),$$
and then the result in \cite{Z} and the isomorphism above imply the contractibility of $D^*(z)$.\\
\indent
For $g$ defined in \eqref{gFantappie} we have
\begin{equation*}
P_j\left(\frac{\partial}{\partial\xi}\right)(g)
=(-1)^{\deg{P}_j-1}\left(\deg{P}_j-1\right)!\mu\left(\frac{P_j(z)}{\langle\xi\cdot z\rangle^{\deg{P}_j}}\right)=0,
\end{equation*}
which concludes the proof of inclusion for ${\cal F}$.\\
\indent
To prove the inclusion for the image of ${\cal R}_V$ we consider for an arbitrary residual current
$\phi\in Z^{(0,n-m-1)}\left(V_D,\omega_V^{\circ}\right)$ the analytic functional on $H(G)$
\begin{equation*}
\mu^{\phi}(h)=\sum_{\alpha=1}^N\lim_{t\to 0}
\int_{bD_{\delta}\cap T^{\epsilon}_{\left\{{\bf F}^{(\alpha)}\right\}}(t)}
\vartheta_{\alpha}(z)\frac{h(z)\cdot\Phi_{\alpha}(z)}
{\prod_{k=1}^m F^{(\alpha)}_k(z)}.
\end{equation*}
From Proposition~\ref{1.7.6(2)} we obtain that $\mu^{\phi}(h)=0$ for any $h \in H^0(G,{\cal I})$,
and therefore $\mu^{\phi}$ defines a functional on $H^0(G, {\cal O}/{\cal I})$. From equality
\eqref{RadonFantappie} we obtain equality
$${\cal F}_V[\mu^{\phi}]=\left(2\pi i\right)^{m+1}{\cal R}_V[\phi],$$
which implies the inclusion
$$\text{Image}\ {\cal R}_V\subseteq\text{Image}\ {\cal F}_V$$
and equality \eqref{R-VEquality}.\\
\indent
To conclude the proof of inclusion for the image of ${\cal R}_V$ we have to prove equality
\begin{equation}\label{muZero}
\mu^{\phi}(h)=0
\end{equation}
for an arbitrary $h\in H^0\left(\C P^n,{\cal O}/{\cal I}\right)$. To prove this equality we assume that
in every $U_{\alpha}$ function $h$ is defined in some neighborhood of $V\cap U_{\alpha}$ and consider
a sequence of nonnegative functions $\eta_{\nu}\in {\cal E}_c(D)$ approximating the characteristic
function of $D_{\delta}$ as $\nu\to 0$. Then, applying the Stokes' formula in each $U_{\alpha}$
we obtain the equality
\begin{multline*}
\sum_{\alpha=1}^N\lim_{t\to 0}\int_{D\cap T^{\epsilon}_{\left\{{\bf F}^{(\alpha)}\right\}}(t)}
\vartheta_{\alpha}(z)\bar\partial\eta_{\nu}(z)\wedge\frac{h(z)\cdot\Phi_{\alpha}(z)}
{\prod_{k=1}^m F^{(\alpha)}_k(z)}\\
+\sum_{\alpha=1}^N\lim_{t\to 0}\int_{D\cap T^{\epsilon}_{\left\{{\bf F}^{(\alpha)}\right\}}(t)}
\eta_{\nu}(z)\bar\partial\vartheta_{\alpha}(z)\wedge\frac{h(z)\cdot\Phi_{\alpha}(z)}
{\prod_{k=1}^m F^{(\alpha)}_k(z)}=0,
\end{multline*}
which is transformed into equality
\begin{equation*}
\sum_{\alpha=1}^N\lim_{t\to 0}\int_{D\cap T^{\epsilon}_{\left\{{\bf F}^{(\alpha)}\right\}}(t)}
\vartheta_{\alpha}(z)\bar\partial\eta_{\nu}(z)\wedge\frac{h(z)\cdot\Phi_{\alpha}(z)}
{\prod_{k=1}^m F^{(\alpha)}_k(z)}=0
\end{equation*}
after application of Lemma~\ref{ZeroIntegral}.\\
\indent
Passing to the limit as $\nu\to 0$ in the equality above we obtain equality \eqref{muZero}.
\qed
%***************************** end of Proposition F-RInclusion **************************

\section{Kernels of ${\cal R}_V$ and ${\cal F}_V$.}\label{KernelDescription}

\indent
In this section we describe the kernels of  ${\cal F}_V$ and ${\cal R}_V$. In the next proposition we prove the triviality
of the kernel of ${\cal F}_V$.

%********************************************************************************
%******* Proposition F-Kernel  ***************************************************
\begin{proposition}\label{F-Kernel}
For the Fantappi\'e transform defined in \eqref{Fantappie} we have
\begin{equation}\label{F-Kernel=0}
\text{Ker}\ {\cal F}_V=\{0\}.
\end{equation}
\end{proposition}
%***************************************************************************************
\indent
{\bf Proof.} To prove property \eqref{F-Kernel=0} we use the linear concavity of $D$ and contractibility of $D^*(z)$
for every $z\in D$, and obtain as in Proposition~\ref{F-RInclusion} the connectedness of $D^*$.
Then using the connectedness of $D^*$ and the Cauchy-Fantappi\'e-Leray integral formula on $G$ (see \cite{L}) we obtain the density of the set of functions
$$\left\{\frac{1}{\xi_0+\sum_{j=1}^n\xi_ju_j}\right\}_{\xi\in D^*}$$
in $H(G)$, where  we used the assumption $D\supset\{z_0=0\}$ and changed variables in $G$ to $u_j=z_j/z_0$.\\
\indent
Then from equality ${\cal F}_V\left[\mu\right](z_0/\langle \xi\cdot z\rangle)=0$ we obtain the equality
$\mu=0$.
\qed\\
%******* End of Proposition F-Kernel  ***************************************************

\indent
In the proposition below we prove the necessity of the condition on $\text{Ker}\ {\cal R}_V$
in the statement (i) of Theorem~\ref{RadonAction}.

%********************************************************************************
%******* Proposition KernelNecessity  ***************************************************
\begin{proposition}\label{KernelNecessity}
Let $V\subset\C P^n$ be a locally complete intersection subvariety, let $D\subset\C P^n$
be a linearly concave domain. If a residual current $\phi\in Z^{(0,n-m-1)}\left(V_D,\omega_V^{\circ}\right)$ 
admits an extension to $\C P^n$ as a $\bar\partial$-closed residual current supported on $V$, then ${\cal R}_V[\phi]=0$.
\end{proposition}
%***************************************************************************************
\indent
{\bf Proof.} Let $\phi\in Z^{(0,n-m-1)}\left(V_D,\omega_V^{\circ}\right)$ be the restriction of
a $\bar\partial$-closed residual current on $V$. Then from equality \eqref{BoundaryEquality}
in Lemma~\ref{TubeIndependence} we obtain
\begin{equation*}
{\cal R}_V[\phi](\xi)=\frac{1}{(2\pi i)^{m+1}}\sum_{j=0}^n\left(\sum_{\alpha=1}^N\lim_{t\to 0}
\int_{bD_{\delta}\cap T^{\epsilon}_{\left\{{\bf F}^{(\alpha)}\right\}}(t)}
\vartheta_{\alpha}(z)\frac{z_j\cdot\Phi_{\alpha}(z)}
{\langle\xi\cdot z\rangle\cdot\prod_{k=1}^m F^{(\alpha)}_k(z)}\right)d\xi_j.
\end{equation*}
\indent
We choose an open domain $U_1\subset G$ from the cover
$\cup_{\alpha=1}^N U_{\alpha}$ of $G$ such that
$$U_1=\left\{z\in G:\ \tau(z)<0\right\}$$
for a function $\tau\in{\cal E}(G)$. Then we consider for a fixed $\mu>0$ a family of smooth
nonnegative functions $\eta_{\nu}$ with compact support such that
$$\eta_{\nu}(z)=\begin{cases}
0\ \mbox{if}\ \tau(z)<-\mu-\nu,\ \mbox{or}\ \rho_{\delta}(z)<-\nu\\
1\ \mbox{if}\ \tau(z)>-\mu+\nu, \mbox{and}\ \rho_{\delta}(z)>\nu.
\end{cases}$$
\indent
As in Lemma~\ref{ZeroIntegral} we have the equality
$$\sum_{\alpha=1}^N\lim_{t\to 0}\int_{T^{\epsilon}_{\left\{{\bf F}^{(\alpha)}\right\}}(t)}
\vartheta_{\alpha}(z)\bar\partial\eta_{\nu}(z)\wedge\frac{z_j\cdot\Phi_{\alpha}(z)}
{\langle\xi\cdot z\rangle\cdot\prod_{k=1}^m F^{(\alpha)}_k(z)}=0,$$
which, after passing to the limit as $\nu\to 0$ implies the equality
\begin{multline*}
\sum_{\alpha=1}^N\lim_{t\to 0}
\int_{bD_{\delta}\cap T^{\epsilon}_{\left\{{\bf F}^{(\alpha)}\right\}}(t)}
\vartheta_{\alpha}(z)\frac{z_j\cdot\Phi_{\alpha}(z)}
{\langle\xi\cdot z\rangle\cdot\prod_{k=1}^m F^{(\alpha)}_k(z)}\\
=\sum_{\alpha=1}^N\lim_{t\to 0}
\int_{\left\{z\in U_1:\ \tau(z)=-\mu\right\}\cap
T^{\epsilon}_{\left\{{\bf F}^{(\alpha)}\right\}}(t)}
\vartheta_{\alpha}(z)\frac{z_j\cdot\Phi_{\alpha}(z)}
{\langle\xi\cdot z\rangle\cdot\prod_{k=1}^m F^{(\alpha)}_k(z)}.
\end{multline*}
\indent
Choosing the partition of unity such that $\vartheta_1\Big|_{\left\{z\in U_1:\ \tau(z)\leq-\mu\right\}}\equiv 1$ we obtain
\begin{multline}\label{R-CompleteBoundary}
{\cal R}_V[\phi](\xi)=\frac{1}{(2\pi i)^{m+1}}\sum_{j=0}^n\left(\sum_{\alpha=1}^N\lim_{t\to 0}
\int_{bD_{\delta}\cap T^{\epsilon}_{\left\{{\bf F}^{(\alpha)}\right\}}(t)}\vartheta_{\alpha}(z)
\frac{z_j\cdot\Phi_{\alpha}(z)}{\langle\xi\cdot z\rangle\cdot\prod_{k=1}^m F^{(\alpha)}_k(z)}\right)d\xi_j\\
=\frac{1}{(2\pi i)^{m+1}}\sum_{j=0}^n\left(
\lim_{t\to 0}\int_{\left\{z\in U_1:\ \tau(z)=-\mu\right\}\cap
T^{\epsilon}_{\left\{{\bf F}^{(1)}\right\}}(t)}
\frac{z_j\cdot\Phi_1(z)}
{\langle\xi\cdot z\rangle\cdot\prod_{k=1}^m F^{(1)}_k(z)}\right)d\xi_j.
\end{multline}
\indent
Then applying the Stokes' formula to the form
$$\frac{z_j\cdot\Phi_1(z)}
{\langle\xi\cdot z\rangle\cdot\prod_{k=1}^m F^{(1)}_k(z)}$$
on the manifold
$$\left\{z\in U_1:\ \tau(z)<-\mu\right\}
\cap T^{\epsilon}_{\left\{{\bf F}^{(1)}\right\}}(t)$$
with the boundary
$$\left\{z\in U_1:\ \tau(z)=-\mu\right\}
\cap T^{\epsilon}_{\left\{{\bf F}^{(1)}\right\}}(t),$$
and using Proposition~\ref{1.7.6(2)} we obtain ${\cal R}_V[\phi]=0$.
\qed\\

{\bf Remark.} The referee has drawn our attention to the fact that Proposition~\ref{KernelNecessity}
must be valid for any current $\phi\in Z^{(0,n-m-1)}\left(V_D,\omega_V^{\circ}\right)$ admitting an extension to $\C P^n$
as a $\bar\partial$-closed current. This is indeed true and can be reduced to the following statement:\\
\indent
{\it
If a current $\phi\in Z^{(0,n-m-1)}\left(V_D,\omega_V^{\circ}\right)\subset
\Gamma\left(D,{\cal K}^{(n,n-1)}\right)$ is $\bar\partial$-cohomologically equivalent to a
$\bar\partial$-closed form $\Phi\in C^{(n,n-1)}(D)$,
then ${\cal R}_V[\phi]={\cal R}[\Phi]$, where ${\cal R}[\Phi]$ is the standard Radon transform of $\Phi$ defined
using the manifold of incidence
$$\Big\{(\xi,z)\in D^*\times D: \langle \xi\cdot z\rangle=0\Big\}.$$
(See similar statement for a reduced $V$ on p.242 in \cite{He1}.)
}%it

\indent
In the next Proposition we prove the sufficiency of the condition in the statement (i) of Theorem~\ref{RadonAction}.

%***************************************************************************************
%******* Proposition KernelSufficiency  *******************************************************
\begin{proposition}\label{KernelSufficiency}\ If a $\bar\partial$-closed residual current $\phi$ on $V_D$ satisfies ${\cal R}_V[\phi]=0$, then $\phi$ is the restriction to $V_D$ of a $\bar\partial$-closed residual current on $V$.
\end{proposition}
%***************************************************************************************
\indent
{\bf Proof.} Without loss of generality we may assume that the $\bar\partial$-closed forms $\Phi_{\alpha}$
associated with $\phi$ are defined in some linearly concave domain $D_{-\delta}\supset D$. 
We fix $\nu$ such that $\delta>\nu>0$ and extend current $\phi$ into $\C P^n$ by extending the forms
$\Phi_{\alpha}$ by the formula $\vartheta\Phi_{\alpha}$, where
$$\vartheta(z)=\begin{cases}
1\ \mbox{if}\ z\in D_{-\nu},\vspace{0.1in}\\
0\ \mbox{if}\ z\notin D_{-\delta},
\end{cases}$$
is a smooth function. Then we consider current $\psi$ defined on $\mathring{G}$ by the formula
\begin{equation}\label{psiCurrent}
\psi\left(f\right)=\sum_{\alpha=1}^N\lim_{t\to 0}
\int_{T^{\epsilon}_{\left\{{\bf F}^{(\alpha)}\right\}}(t)}
\vartheta_{\alpha}(z)f(z)\frac{\bar\partial\vartheta(z)\wedge\Phi_{\alpha}(z)}{\prod_{k=1}^m F^{(\alpha)}_k(z)}
\end{equation}
for $f\in {\cal E}_c(\mathring{G})$.\\
\indent
Using the Stokes' formula, Lemma~\ref{ZeroIntegral}, and Proposition~\ref{1.7.6(2)} we obtain the following equality
\begin{multline*}
\sum_{\alpha=1}^N\lim_{t\to 0}\int_{T^{\epsilon}_{\left\{{\bf F}^{(\alpha)}\right\}}(t)}
\vartheta_{\alpha}(z)f(z)\frac{\bar\partial\vartheta(z)\wedge\Phi_{\alpha}(z)}{\prod_{k=1}^m F^{(\alpha)}_k(z)}\\
=-\sum_{\alpha=1}^N\lim_{t\to 0}\int_{T^{\epsilon}_{\left\{{\bf F}^{(\alpha)}\right\}}(t)}
\vartheta_{\alpha}(z)\bar\partial f(z)\wedge\frac{\vartheta(z)\Phi_{\alpha}(z)}{\prod_{k=1}^m F^{(\alpha)}_k(z)}\\
-\sum_{\alpha=1}^N\lim_{t\to 0}\int_{T^{\epsilon}_{\left\{{\bf F}^{(\alpha)}\right\}}(t)}
 f(z)\vartheta(z)\bar\partial\vartheta_{\alpha}(z)\wedge\frac{\Phi_{\alpha}(z)}{\prod_{k=1}^m F^{(\alpha)}_k(z)}\\
-\sum_{\alpha=1}^N\lim_{t\to 0}\int_{T^{\epsilon}_{\left\{{\bf F}^{(\alpha)}\right\}}(t)}
 f(z)\vartheta(z)\vartheta_{\alpha}(z)\wedge\frac{\bar\partial\Phi_{\alpha}(z)}{\prod_{k=1}^m F^{(\alpha)}_k(z)}
\end{multline*}
\begin{equation*}
=-\sum_{\alpha=1}^N\lim_{t\to 0}\int_{T^{\epsilon}_{\left\{{\bf F}^{(\alpha)}\right\}}(t)}
\vartheta_{\alpha}(z)\bar\partial f(z)\wedge\frac{\vartheta(z)\Phi_{\alpha}(z)}{\prod_{k=1}^m F^{(\alpha)}_k(z)},
\end{equation*}
i.e. $\psi$ is a current with compact support in $\mathring{G}$ satisfying the condition
$$\psi=\bar\partial\left(\vartheta\phi\right).$$
\indent
Considering the extension of $\psi$ to the space of holomorphic functions on $\mathring{G}$ and using the
Stokes' formula we obtain 
\begin{equation}\label{AnalyticOrthogonality}
\sum_{\alpha=1}^N\lim_{t\to 0}\int_{T^{\epsilon}_{\left\{{\bf F}^{(\alpha)}\right\}}(t)}
\vartheta_{\alpha}(z)h(z)\frac{\bar\partial\vartheta(z)\wedge\Phi_{\alpha}(z)}{\prod_{k=1}^m F^{(\alpha)}_k(z)}
=\sum_{\alpha=1}^N\lim_{t\to 0}\int_{bD_{-\nu}\cap T^{\epsilon}_{\left\{{\bf F}^{(\alpha)}\right\}}(t)}
\vartheta_{\alpha}(z)h(z)\frac{\Phi_{\alpha}(z)}{\prod_{k=1}^m F^{(\alpha)}_k(z)}
\end{equation}
for a holomorphic $h\in H(\mathring{G})$.\\
\indent
Using condition ${\cal R}_V[\phi]=0$ and introducing variables
$$u_j=\frac{z_j}{z_0}\hspace{0.1in}\mbox{for}\ j=1,\dots,n,$$
in the neighborhood $\left\{z_0\neq 0\right\}$ we obtain the equality
$${\cal R}_V[\phi]_0(\xi)=\frac{1}{(2\pi i)^{m+1}}\sum_{\alpha=1}^N
\lim_{t\to 0}\int_{bD_{-\nu}\cap T^{\epsilon}_{\left\{{\bf F}^{(\alpha)}\right\}}(t)}
\vartheta_{\alpha}(u)\frac{\Phi_{\alpha}(u)}{\left(\xi_0+\sum_{l=1}^n\xi_l\cdot u_l\right)
\cdot\prod_{k=1}^m F^{(\alpha)}_k(u)}=0$$
for arbitrary $\xi\in D^*$.\\
\indent
From the linear concavity of $D$ and contractibility of $D^*(z)$ for every $z\in D$, which we pointed out
above in Proposition~\ref{F-Kernel}, we obtain the connectedness of $D^*$.
Then again using the connectedness of $D^*$ and the Cauchy-Fantappi\'e-Leray integral formula
on $G$ (see \cite{L}) we obtain the density of the set of functions
$$\left\{\frac{1}{\xi_0+\sum_{j=1}^n\xi_ju_j}\right\}_{\xi\in D^*_{\delta}}$$
in $H(\mathring{G})$. Then the equality
$$\sum_{\alpha=1}^N\lim_{t\to 0}\int_{bD_{-\nu}\cap T^{\epsilon}_{\left\{{\bf F}^{(\alpha)}\right\}}(t)}
\vartheta_{\alpha}(u)\frac{h(u)\Phi_{\alpha}(u)}{\prod_{k=1}^m F^{(\alpha)}_k(u)}=0$$
holds for an arbitrary $h\in H(\mathring{G})$, which implies, according to \eqref{AnalyticOrthogonality},
the equality
\begin{equation}\label{psiZero}
\psi\left(h\right)=\sum_{\alpha=1}^N\lim_{t\to 0}\int_{T^{\epsilon}_{\left\{{\bf F}^{(\alpha)}\right\}}(t)}
\vartheta_{\alpha}(z)h(z)\frac{\bar\partial\vartheta(z)\wedge\Phi_{\alpha}(z)}{\prod_{k=1}^m F^{(\alpha)}_k(z)}=0.
\end{equation}
\indent
From the Serre-Malgrange duality (see \cite{Mal2}, \cite{S}) one can obtain (see  \cite{DGSY}, \S 2,  Lemma 2.2) that
\begin{equation}\label{FactorEquality}
H^0\left(\mathring{G},{\cal O}/{\cal I}\right)^{\prime}
=\Gamma_c\left(\mathring{G},{\cal K}_{\cal I}^{n,n}\right)/
\bar\partial \left\{\Gamma_c\left(\mathring{G},{\cal K}_{\cal I}^{n,n-1}\right)\right\},
\end{equation}
where ${\cal I}$ is the sheaf of ideals defined by the polynomials $\left\{P_1,\dots,P_r\right\}$ and
${\cal K}_{\cal I}^{p,q}$ is the sheaf of germs of currents $\gamma^{(p,q)}$ on $\mathring{G}$
with compact support in $V$ such that for any open subset $U\subset \mathring{G}$
the current $\gamma$ satisfies
$$\gamma\left(g\cdot f\right)=0$$
for any $g\in H^0\left(U, {\cal I}\right)$ and $f\in {\cal E}_c^{(n-p,n-q)}\left(U\right)$.\\
\indent
From equality \eqref{FactorEquality} applied to the current $\psi$ defined in \eqref{psiCurrent}
using \eqref{psiZero}, we obtain the existence of
$\beta\in \Gamma_c\left(\mathring{G},{\cal K}_{\cal I}^{n,n-1}\right)$ satisfying
$$\bar\partial\beta=\psi,$$
and therefore, the current $\beta-\vartheta\phi$ is an extension of the current $\phi$ into $G$
as a $\bar\partial$-closed current. The existence of such current is precisely the appropriate
modification of the statement of Theorem~\ref{RadonAction} mentioned in the remark to this theorem.
Namely, if $m<n-1$, $V\subset D$ is a locally complete intersection in $D$, and a $\bar\partial$-closed
residual current $\phi$ on $V_D$ satisfies ${\cal R}_V[\phi]=0$, then $\phi$ admits a $\bar\partial$-closed
extension to a current $\gamma$ on $\C P^n$ satisfying
$$\gamma\left(g\cdot f\right)=0$$
for any $g\in H^0\left(U, {\cal I}\right)$ and $f\in {\cal E}_c^{(n-p,n-q)}\left(U\right)$.\\
\indent
If $V$ is a locally complete intersection in $\C P^n$, then a residual current extension can be found. In this case
using the partition of unity $\left\{\vartheta_{\alpha}\right\}_{\alpha=1}^N$ we rewrite the last equality as
$$\sum_{\alpha=1}^N\bar\partial\left(\vartheta_{\alpha}\beta\right)=\psi$$
with currents $\vartheta_{\alpha}\beta$ having compact supports in $U_{\alpha}$ and satisfying
$$\vartheta_{\alpha}\beta\left(g\cdot f^{(0,1)}\right)=0$$
for any $g\in H^0\left(U, {\cal I}\right)$ and $f^{(0,1)}\in {\cal E}_c^{(0,1)}\left(U\right)$.\\
\indent
Using then the result of Dickenstein-Sessa \cite{DS} motivated by Palamodov \cite{P2}
(see also Theorem 3.4 from \cite{DGSY}) we obtain the existence in $\left\{U_{\alpha}\right\}_{\alpha=1}^N$
of a collection of residual currents $\theta_{\alpha}$ with compact support in $U_{\alpha}$ of the form
$$\theta_{\alpha}(f)=\lim_{t\to 0}\int_{T^{\epsilon}_{\left\{{\bf F}^{(\alpha)}\right\}}(t)}
\frac{f(u)\wedge \Theta_{\alpha}(u)}{\prod_{k=1}^m F^{(\alpha)}_k(u)},$$
where $\Theta_{\alpha}$ are $\bar\partial$-closed forms of type $(n,n-m-1)$ in some neighborhood
of $U_{\alpha}\cap V$ with compact support in $U_{\alpha}$, such that
$$\bar\partial\left(\vartheta_{\alpha}\beta-\theta_{\alpha}\right)=0.$$
\indent
Therefore, the current $\vartheta\phi-\theta$ is an extension of current $\phi$ into $G$
as a $\bar\partial$-closed residual current.\qed
%*******End Proposition KernelSufficiency  ***************************************************

\section{Images of ${\cal F}_V$ and ${\cal R}_V$.}\label{F-RImages}
\indent
In this section we complete the proof of Theorem~\ref{RadonAction} by proving
the second part of statement (ii), namely the inclusions
\begin{multline}\label{RightInclusionProperty}
\text{Image}\ {\cal F}_V
\supseteq \left\{f\in H^{(1,0)}(D^*): f=dg\ \mbox{with}\ g\in H^0(D^*)\ \mbox{such that}\
\left\{P_k\left(\frac{\partial}{\partial\xi}\right)g=0\right\}_1^r\right\},\\
\text{Image}\ {\cal R}_V\supseteq\left\{f\in \text{Image}\ {\cal F}_V:\ f={\cal F}_V[\mu],\ \mbox{where}\
\mu(h)=0\ \text{for}\ \forall  h\in H^0\left(\C P^n,{\cal O}/{\cal I}\right)\right\},
\end{multline}
and statement (iii) of this theorem.

\indent
In the proposition below we prove the inclusion above for the image of the Fantappi\'e transform.

%********************************************************************************
%******* Proposition FantappieMartineau  ************************************************
\begin{proposition}\label{FantappieMartineau}
Under the hypotheses of Theorem~\ref{RadonAction} for any
$f=dg\in H^{(1,0)}\left(D^*\right)$ with $g$ satisfying equations
\begin{equation}\label{System}
P_1\left(\frac{\partial}{\partial\xi}\right)g=\cdots=P_r\left(\frac{\partial}{\partial\xi}\right)g=0,
\end{equation}
there exists a linear functional $\mu\in H^0(G,{\cal O}/{\cal I})^{\prime}$,
such that ${\cal F}_V[\mu]=f$.
\end{proposition}
%********************************************************************************
\indent
{\bf Proof.} To prove the proposition we use the following version of the
Martineau's (see \cite{Mar2}) inversion formula from \cite{GH}.

%********************************************************************************
%******* Proposition DLerayMartineau  ************************************************
\begin{proposition}\label{DLerayMartineau}(Generalized Martineau inversion formula.
\cite{Mar2}, \cite{GH}.)\ Let $D\subset \C P^n$ be a linearly concave domain such that
$D^*\subset \left\{\xi_0\neq 0\right\}$, and let $g\in H^0(D^*)$ be
a holomorphic function of homogeneity $0$ on $D^*$.\\
\indent
Let  $\mu^g$ be the linear functional on $H(G)$ defined by the formula (see \cite{Mar2,GH})
\begin{equation}\label{mu_g}
\mu^g(h)=\int_{bG_{\nu}}h\cdot \Omega_g,
\end{equation}
where
$$\Omega_g(z)=\frac{(-1)}{(2\pi i)^n}\ \frac{\partial^n g}{\partial\xi_0^n}(\eta(z))
\omega^{\prime}\left(\eta(z)\right)\bigwedge_{j=1}^nd\left(\frac{z_j}{z_0}\right),$$
$$\omega^{\prime}\left(\eta\right)=\sum_{j=1}^n(-1)^j\eta_jd\eta_1\wedge
\stackrel{\stackrel{j}{\vee}}{\cdots}\wedge d\eta_n,$$
and a map $\eta:bG_{\nu}\to D^*$ satisfies
$\langle\eta(z)\cdot z\rangle=0$ for $z\in bG_{\nu}$.\\
\indent
Then the following equality holds:
\begin{equation}\label{DMartineauFormula}
{\cal F}[\mu^g](\xi)=dg(\xi),
\end{equation}
or
$$\frac{(-1)}{(2\pi i)^n}\int_{bG_{\nu}}\frac{z_k}{\langle\xi\cdot z\rangle}
\frac{\partial^{n}g}{\partial\xi_0^{n}}(\eta(z))
\omega^{\prime}\left(\eta(z)\right)\bigwedge_{j=1}^nd\left(\frac{z_j}{z_0}\right)
=\frac{\partial g}{\partial\xi_k}(\xi)\ \ \mbox{for}\ k=0,\dots,n$$
for $\xi\in D^*$.
\end{proposition}
\qed\\
%********************************************************************************
\indent
Using Proposition~\ref{DLerayMartineau} we construct for an arbitrary $g\in H^0(D^*)$ the current $\mu^g$
satisfying equality \eqref{DMartineauFormula}. To prove that $\mu^g(h)=0$ for any $h \in H^0(G,{\cal I})$,
and that therefore $\mu^g$ defines a functional on $H^0(G, {\cal O}/{\cal I})$
we use the assumption on $g$, to obtain the equality
$$(-1)^{1+\deg{P}_k}\left(\deg{P}_k\right)!\cdot
\int_{bG_{\nu}}\frac{z_0}{\langle\xi\cdot z\rangle^{1+\deg{P}_k}}P_k(z)\Omega_g(z)
=P_k\left(\frac{\partial}{\partial\xi}\right)
\left[\frac{\partial g}{\partial\xi_0}\right]=0.$$
Then from the connectedness of $D^*$ (see discussion in Proposition~\ref{F-Kernel}), and therefore
the density of the set of functions
$$\left\{\frac{z_0}{\langle\xi\cdot z\rangle^{1+\deg{P}_k}}\right\}_{\xi\in D^*}$$
in the space $H^0(G)$, we obtain the equality
\begin{equation}\label{ZeroonH}
\mu^g(h\cdot P_k)=\frac{(-1)}{(2\pi i)^n}\int_{bG_{\nu}}h(z)\cdot P_k(z)
\cdot\frac{\partial^n g}{\partial\xi_0^n}(\eta(z))
\omega^{\prime}\left(\eta(z)\right)\bigwedge d\left(\frac{z_j}{z_0}\right)=0
\end{equation}
for an arbitrary $h\in H^0(G)$.\\
\qed
%******* End Proposition FantappieMartineau  *********************************************

We prove the second inclusion from \eqref{RightInclusionProperty}
and statement (iii) of Theorem~\ref{RadonAction} in the following proposition.

%************************************************************************************
%******* Proposition RightInclusion  *******************************************************
\begin{proposition}\label{RightInclusion} Under the hypotheses of Theorem~\ref{RadonAction} for any
$f=dg\in H^{(1,0)}\left(D^*\right)$ with $g$ satisfying equations \eqref{System}
and $\mu^g$ constructed in Proposition~\ref{FantappieMartineau} satisfying
\begin{equation}\label{RadonCondition}
{\cal F}[\mu^g]=dg,\ \text{and}\ \mu^g(h)=0\ \text{for}\ \forall  h\in H^0\left(\C P^n,{\cal O}/{\cal I}\right),
\end{equation}
there exists a residual current $\phi\in Z^{(0,n-m-1)}\left(V_D,\omega_V^{\circ}\right)$,
such that ${\cal R}_V[\phi]=f$.\\
\indent
Such current also exists if $V$ is connected in the sense that $\dim H^0(V,{\cal O}/{\cal I})=1$.
\end{proposition}
%**************************************************************************************
\indent
{\bf Proof.} To construct a $\bar\partial$-closed residual current with support on $V_D$,
such that its Radon transform coincides with $dg$, we need an identification described below.\\
\indent
First we consider the following equality of Hartshorne (see \cite{Ha}, Ch.III, Corollary 7.7, Theorem 7.11),
specifying the results of Serre \cite{S}, Grothendieck \cite{Gro}, Ramis, Ruget, Verdier \cite{RR}, \cite{RRV}
for locally complete intersections
\begin{equation}\label{DualizingEquality}
H^0\left(V,{\cal O}/{\cal I}\right)^{\prime}\cong H^{n-m}\left(V, \omega_V^{\circ}\right),
\end{equation}
where ${\cal I}$ is the sheaf of germs of ideals corresponding to $V$,
and $\omega_V^{\circ}=\omega_{\C P^n}\otimes \det N(V)^{-1}$ is the {\it dualizing sheaf} of $V$
defined earlier in \eqref{Dualizing}.\\
\indent
Using the exactness of the $\bar\partial$-complex of sheaves
\begin{equation*}
0\to {\cal O}/{\cal I}\otimes\omega_{\C P^n}\to {\cal O}/{\cal I}\otimes{\cal E}^{(n,0)}
\stackrel{\bar\partial}{\to}\cdots\stackrel{\bar\partial}{\to} {\cal O}/{\cal I}\otimes{\cal E}^{(n,n)}\to 0,
\end{equation*}
which follows from the Malgrange's theorem on ${\cal O}$-flatness of ${\cal E}$
(see \cite{Mal1}, $n^{\circ} 25$, Th. 2), we obtain the equality
\begin{multline}\label{Cech-Dolbeault}
H^{n-m}\left(V, \omega_V^{\circ}\right)
\cong H^{n-m}_{\bar\partial}\left(V, \omega_V^{\circ}\right)\\
\cong\frac{\left\{\phi\in {\cal E}^{(n,n-m)}\left(U,\det N(V)^{-1}\right):\
\bar\partial\phi\in {\cal I}\otimes {\cal E}^{(n,n-m+1)}\left(U,\det N(V)^{-1}\right)\right\}}
{\left\{\phi\in \bar\partial {\cal E}^{(n,n-m-1)}\left(U,\det N(V)^{-1}\right)
+{\cal I}\otimes {\cal E}^{(n,n-m)}\left(U,\det N(V)^{-1}\right)\right\}}
\end{multline}
for a small enough neighborhood $U\supset V$.\\
\indent
On the other hand, for any representative
$\Phi\in H^{n-m}_{\bar\partial}\left(V, \omega_V^{\circ}\right)$ using the Coleff-Herrera theory
we can construct a linear functional on $H^0\left(V,{\cal O}/{\cal I}\right)$ by the formula
\begin{equation}\label{Herrera}
\langle\phi,h\rangle=\sum_{\alpha=1}^N\lim_{t\to 0}\int_{T^{\epsilon}_{\left\{{\bf F}^{(\alpha)}\right\}}(t)}
\vartheta_{\alpha}(z)\frac{h(z)\Phi_{\alpha}(z)}{\prod_{k=1}^m F^{(\alpha)}_k(z)},
\end{equation}
explicitly defining the isomorphism in \eqref{DualizingEquality}.\\
\indent
Continuing then with the construction of the sought current we observe that for an arbitrary fixed
$\delta>0$ and the analytic functional $\mu^g$ on $H(\mathring{G}_{\delta})$
defined in \eqref{mu_g} we can use equality \eqref{FactorEquality} and obtain the existence
of a current $\psi^{(\delta)}\in\Gamma_c\left(\mathring{G}_{\delta},{\cal K}_{\cal I}^{n,n}\right)$
with support in $V\cap \mathring{G}_{\delta}$, coinciding with the analytic functional $\mu^g$
on $H(\mathring{G}_{\delta})$ defined in \eqref{mu_g}.
Considering current $\psi^{(\delta)}$ as a current on $V$ and using equality \eqref{ZeroonH} we obtain the
existence of a $\bar\partial$-closed differential form
$$\Psi^{(\delta)}\in {\cal E}^{(0,n-m)}\left(V, \omega_V^{\circ}\right)$$
corresponding to $\mu^g$  by equality \eqref{DualizingEquality} and such that
\begin{equation}\label{mu-g-current}
\psi^{(\delta)}(h)=\sum_{\alpha=1}^N\lim_{t\to 0}\int_{T^{\epsilon}_{\left\{{\bf F}^{(\alpha)}\right\}}(t)}
\vartheta_{\alpha}(z)\frac{h(z)\Psi^{(\delta)}_{\alpha}(z)}{\prod_{k=1}^m F^{(\alpha)}_k(z)}
\end{equation}
by equality \eqref{Herrera}.\\
\indent
Using condition \eqref{RadonCondition} for functional $\mu^g$ and equality \eqref{Cech-Dolbeault}
we obtain that the functional in \eqref{mu-g-current} is equal to zero, i.e.
$\psi^{(\delta)}=0$ in $H^{n-m}_{\bar\partial}\left(V, \omega_V^{\circ}\right)$. Therefore, there exists an element
$$\Theta^{(\delta)}\in {\cal E}^{(0,n-m-1)}\left(U, \omega_V^{\circ}\right)$$
in some neighborhood $U$ of $V$ such that
$$\bar\partial \Theta^{(\delta)}\Big|_V=\Psi^{(\delta)}\Big|_V.$$
\indent
Since $\Psi^{(\delta)}$ has a support in $G_{\delta}$, it follows that the restriction of the form $\Theta^{(\delta)}$
to $D_{\delta}$ is a $\bar\partial$-closed form on $V\cap D_{\delta}$, and the current
$$\theta^{(\delta)}\left(\gamma^{(0,1)}\right)=\sum_{\alpha=1}^N\lim_{t\to 0}
\int_{T^{\epsilon}_{\left\{{\bf F}^{(\alpha)}\right\}}(t)}\vartheta_{\alpha}(\zeta)
\frac{\gamma\wedge\Theta^{(\delta)}_{\alpha}(\zeta)}{\prod_{k=1}^m F^{(\alpha)}_k(\zeta)}$$
is a $\bar\partial$-closed closed residual current in $D_{\delta}$ with support in $V\cap D_{\delta}$.\\
\indent
Applying the Radon transform to the current $\theta^{(\delta)}$ and using equality \eqref{BoundaryEquality}
we obtain the equality
\begin{multline}\label{DeltaRadon}
{\cal R}_V\left[(2\pi i)^{m+1}\cdot\theta^{(\delta)}\right]=\sum_{j=0}^n\left(\sum_{\alpha=1}^N
\lim_{t\to 0}\int_{bD_{\delta}\cap T^{\epsilon}_{\left\{{\bf F^{(\alpha)}}\right\}}(t)}\vartheta_{\alpha}(z)\frac{z_j\cdot\Theta^{(\delta)}_{\alpha}(z)}
{\langle\xi\cdot z\rangle\cdot\prod_{k=1}^m F^{(\alpha)}_k(z)}\right)d\xi_j\\
=\sum_{j=0}^n\left(\sum_{\alpha=1}^N
\lim_{t\to 0}\int_{G\cap T^{\epsilon}_{\left\{{\bf F^{(\alpha)}}\right\}}(t)}\frac{z_j\cdot\Psi^{(\delta)}_{\alpha}(z)}
{\langle\xi\cdot z\rangle\cdot\prod_{k=1}^m F^{(\alpha)}_k(z)}\right)d\xi_j\\
=\sum_{j=0}^n\psi^{(\delta)}\left(\frac{z_j}{\langle\xi\cdot z\rangle}\right)d\xi_j
=\sum_{j=0}^n\mu^g\left(\frac{z_j}{\langle\xi\cdot z\rangle}\right)d\xi_j
={\cal F}[\mu^g](\xi)=dg(\xi).
\end{multline}
\indent
Using the same arguments as above we construct currents $\psi^{(\delta^{\prime})}$
and $\theta^{(\delta^{\prime})}$ for an arbitrary $\delta^{\prime}<\delta$. Then, from \eqref{DeltaRadon}
we obtain the equality
$${\cal R}_V\left[\theta^{(\delta)}-\theta^{(\delta^{\prime})}\right](\xi)=0$$
for $\xi\in  D^*_{\delta}$, and therefore, applying Proposition~\ref{KernelSufficiency} to the current
$\theta^{(\delta)}-\theta^{(\delta^{\prime})}$ on $D_{\delta}$ we obtain the existence of a $\bar\partial$-closed current
$\omega^{(\delta)}$ on $V$, such that
$$\theta^{(\delta)}+\omega^{(\delta)}\Big|_{V\cap D_{\delta}}=\theta^{(\delta^{\prime})},$$
and therefore
\begin{equation}\label{BarEquality}
\bar\partial \theta^{(\delta)}=\bar\partial \theta^{(\delta^{\prime})}=\psi^{(\delta^{\prime})}.
\end{equation}
The equality above shows that the support of $\bar\partial \theta^{(\delta)}$ belongs to $G_{\nu}$
with arbitrary $\nu>0$, i.e. the restriction of the constructed residual current $\theta^{(\delta)}$ to $D$
is a $\bar\partial$-closed current satisfying \eqref{DeltaRadon}.\\
\indent
This completes the proof of the second inclusion in \eqref{RightInclusionProperty}.\\

To prove statement (iii) of Theorem~\ref{RadonAction} we notice that if
$\dim H^0\left(\C P^n,{\cal O}/{\cal I}\right)=1$, then using equality
\begin{equation*}
\mu^g\left(1\right)=\mu^g\left(\frac{z_0}{1\cdot z_0+0\cdot z_2+\cdots 0\cdot z_n}\right)
=\left[{\cal F}\mu^g\right]_0(1,0,\dots,0)
=\frac{\partial g}{\partial\xi_0}(1,0,\dots,0)=0,
\end{equation*}
we obtain that functional $\mu^g$ is equal to zero on $H^0\left(\C P^n,{\cal O}/{\cal I}\right)$,
and therefore using  equality \eqref{Cech-Dolbeault} we obtain
that $\psi^{(\delta)}=0$ in $H^{n-m}_{\bar\partial}\left(V, \omega_V^{\circ}\right)$. The rest of the proof in this case
goes exactly as in the proof above.
\qed
%******************************* end of Proposition RightInclusion **********************
%\newpage

\end{document}